\documentclass{amsart}
\usepackage{amsfonts}
\usepackage{amscd}

\newtheorem{theorem}{Theorem}[section]
\newtheorem{lemma}[theorem]{Lemma}
\newtheorem{Cor}[theorem]{Corollary}
\theoremstyle{definition}
\newtheorem{definition}[theorem]{Definition}
\newtheorem{example}[theorem]{Example}
\newcommand{\at}{\symbol{'100}}

\newcommand{\cent}{\mathcal C}
\newcommand{\sN}{\mathcal N}
\newcommand{\z}{\mathcal Z}
\newcommand{\Sch}{\mbox{\rm Sch}}
\newcommand{\pSch}{\mbox{\rm $p$Sch}}
\newcommand{\oP}{\mbox{$o\!P$}}
\newcommand{\f}{\mathfrak}
\newcommand{\fF}{{\mathfrak F}}

\newcommand{\p}{{\mbox{$[p]$}}}
\newcommand{\scp}{{\mbox{$\scriptstyle [p]$}}}
\newcommand{\pd}{{\mbox{$[p]'$}}}
\newcommand{\scpd}{{\mbox{$\scriptstyle [p]'$}}}
\DeclareMathOperator{\Proj}{Proj}
\DeclareMathOperator{\Cov}{Cov}
\DeclareMathOperator{\soc}{Soc}
\DeclareMathOperator{\im}{im}
\DeclareMathOperator{\Hom}{Hom}
\DeclareMathOperator{\ad}{ad}
\DeclareMathOperator{\Res}{Res}
\DeclareMathOperator{\Und}{Under}
\DeclareMathOperator{\Ord}{Ord}

\DeclareMathOperator{\Skel}{Skel}

\DeclareMathOperator{\Ev}{Ev}
\DeclareMathOperator{\pEv}{\mbox{\rm$p$Ev}}

\DeclareMathOperator{\envd}{Envd}
\DeclareMathOperator{\pLoc}{\mbox{\rm$p$Loc}}
\DeclareMathOperator{\quot}{\textsc{q}\mspace{-2mu}}
\DeclareMathOperator{\sdir}{\textsc{r}\mspace{-2mu}}
\DeclareMathOperator{\frat}{{\textsc{e}}_\Phi\mspace{-3mu}}
\DeclareMathOperator{\prim}{\textsc{p}\mspace{-2mu}}

\newcommand{\pCs}{\mbox{$p\f{C}$}}

\title{Schunck classes of soluble restricted Lie algebras}
\author{Donald W. Barnes}
\address{1 Little Wonga Rd.\\Cremorne NSW 2090\\Australia\\e-mail:
donwb{\at}iprimus.com.au}
\thanks{This work was done while the author was an Honorary Associate of the School of
Mathematics and Statistics, University of Sydney}


\begin{document}

\maketitle

\begin{abstract} I set out the theory of Schunck classes and projectors for  soluble restricted     
Lie algebras and investigate its links to the corresponding theory for ordinary soluble Lie algebras     
over a field $F$ of characteristic $p \ne 0$.
\end{abstract}

\section{Introduction}\label{intro}
Schunck classes and formations were originally defined in the context of finite soluble groups
where  they provide families of intravariant subgroups, the $\f{H}$-projectors.  Development of the
analogous theory for finite dimensional soluble  Lie algebras was started  by Barnes and Gastineau-Hills in
\cite{BGH} and extended by Stitzinger in \cite{Stit1} and \cite{Stit2}. For Lie algebras, the notion of
intravariance has to be defined in terms of derivations rather than automorphisms.  It was shown in
Barnes \cite{Frat}, that the Lie algebra $\f{H}$-projectors are intravariant in this sense.

The theories of Schunck classes of soluble Lie algebras and of restricted Lie algebras over a
field $F$ of characteristic $p$ do not fit together smoothly.  The restrictable Lie algebras do not
form a Schunck class.  Every non-zero Schunck class contains all nilpotent algebras,
but not every nilpotent Lie algebra is restrictable.  The two theories however are linked by
Theorem 6.4 of Barnes \cite{extras}, which I quote here for convenience of reference.

\begin{theorem} \label{thm-prep}  Let $(L, [p])$ be a restricted Lie algebra over the
field $F$ of characteristic $p \ne 0$ and suppose that $z^{[p]} = 0$ for all $z$ in the
centre of $L$.  Let $\fF$ be a saturated formation and suppose $S \ne 0$ is subnormal
in $L$ and $S \in \fF$.  Let $V$ be an irreducible $p$-module of $L$.  Then $V$ is
$S\fF$-hypercentral.
\end{theorem}

Developments in the theory of finite groups have led to changes in the terminology which have not so far
been copied in the Lie algebra literature.   In this paper, I follow the terminology which has
become standard in finite group theory, set out in Doerk and Hawkes \cite{DH}.  I construct a theory of
Schunck classes, saturated formations and projectors within the category of finite-dimensional soluble
restricted Lie algebras over a field $F$ of characteristic $p \ne 0$ and investigate its relationship to the
theory for ordinary soluble Lie algebras over the same field.  We shall see that if a  restricted Lie algebra
$(L, \p)$ is in the Schunck class $\f{H}$ which contains all abelian algebras and if $\pd$ is
another $p$-operation on $L$, then $(L, \pd)$ is also in $\f{H}$.  Thus the theory becomes in effect, a
theory of Schunck classes and saturated formations of restrictable soluble Lie algebras.

The development in Doerk and Hawkes is based on operations on classes of groups.  In whatever
category we are working, by a {\em class} of objects, we mean a collection $\f{X}$ of objects of
the category with the property that if $X \in \f{X}$ and $Y \simeq X$, then $Y \in \f{X}$.  The
foundational issue of operations defined on classes is easily avoided.  Because our algebras are all
finite-dimensional over a given field $F$, we can construct a {\em set} of representative algebras such
that every object of our category is isomorphic to at least one of them.  Our operations could then be
defined in terms of these representatives, at the cost of clumsiness of terminology.   I shall not bother
with this.

The arguments used in Doerk and Hawkes \cite{DH}, with the exception of those referring to conjugacy,
are not specific to groups but are readily adapted to other types of algebras subject to the following
conditions.  Firstly, we require that there be an action of an algebra on itself such that the subalgebras
stable under this action are precisely the kernels of homomorphisms.  For groups, this action is
conjugation.  For Lie algebras, it is the adjoint action.  Secondly, we require a finiteness condition,
some measure of ``size'' of an algebra such that proper subalgebras and proper quotients have smaller
size.  Thirdly, we require the solubility condition that for the action of an algebra $A$ on a chief factor
$B/C$, $B$ is contained in the kernel of that action.  These conditions ensure a covering-avoidance
property: if $M$ is a maximal subalgebra of $A$ and $B/C$ is a chief factor of $A$, then $M$ either
covers $B/C$, that is, $M + C\supseteq B$ or avoids it, that is, $M \cap B \subseteq C$.  A fourth
requirement, related to solubility, is that if $B$ is precisely the kernel of the action of $A$ on the chief
fector $B/C$, then $A/C$ splits over $B/C$.

In Sections \ref{prelim} and \ref{frat}, I set out some preliminary results.  Once we
have these, the theory can be developed much as the theory for finite groups or ordinary Lie algebras. 
This is done in Sections \ref{proj}, \ref{hyperc} and \ref{constr}.  The intravariance of projectors
is established in Section \ref{intrav}.  In Section \ref{compar}, we look at the relationship  
between the two theories.   In Section \ref{env}, we consider $p$-envelopes as another link between
them.

Finite dimensional restricted Lie algebras over a field $F$ of characteristic $p > 0$ and
infinitesimal algebraic group schemes over $F$ of height $\le 1$ are equivalent categories.  (See
Demazure and Gabriel \cite[Chapter II, \S 7, Proposition 4.1]{DG}.)  Working in the context of
group schemes, Voigt \cite{Vgt} has obtained some of the preliminary results given below.  For
most of these, the proofs given here are shorter and more elementary.

\section{Preliminaries}\label{prelim}

All algebras considered will be finite-dimensional  over a field $F$ of characteristic
$p \ne 0$.  A restricted Lie algebra, as well as the usual Lie algebra operations, has a
$p$-operation $\p : L \to L$ satisfying $\ad(x)^p = \ad(x^\scp)$, $(\lambda x)^\scp =
\lambda^p x^\scp$ and $(x + y)^\scp = x^\scp + y^\scp + S(x,y)$ where $S(x,y)$ is a function
whose formal definition and basic properties are set out in Strade and
Farnsteiner \cite[Chapter 2]{SF}.   In the special case where $xy=0$, this simplifies to
$(x+y)^\scp = x^\scp + y^\scp$.  Knowledge of the formal definition is not needed for this
paper.  We use the basic properties of restricted Lie algebras set out in Strade and
Farnsteiner \cite[Chapter 2]{SF}.  Construction of examples will be based on the
following theorem of Jacobson (see Jacobson \cite[Theorem V.11, p.190]{Jac} or Strade and
Farnsteiner \cite[Theorem 2.3, p.71]{SF}), of which we shall make frequent use:

\begin{theorem} \label{JacThm} Let $\{a_1, \dots, a_n\}$ be a basis of the Lie algebra $L$
and let $b_1, \dots, b_n$ be elements of $L$ such that $\ad(a_i)^p = \ad(b_i)$ for all
$i$.  Then there exists one and only one $p$-operation \p on $L$ such that $a_i^\scp =
b_i$ for all $i$. \end{theorem}

A Lie algebra $L$ is called restrictable if, for all $x \in L$, $\ad(x)^p$ is an inner
derivation.  By Jacobson's Theorem, if $L$ is restrictable, then there exists at least
one $p$-operation on $L$.  Further, if $L$ is restrictable and $A$ is an ideal of $L$,
then $L/A$ is restrictable.
We shall from time to time, have need to consider restricted Lie algebras $(L, \p)$ and $(L, \pd)$
with the same underlying Lie algebra $L$.  Note that in this situation, $x^\scp$ and $x^\scpd$ differ by
an element of the centre $\z(L)$.  A subalgebra $U$ of $L$ is called a $\p$-subalgebra if $U^\scp
\subseteq U$.  As the centre $\z(L)$ of $(L, \p)$ and so all terms of the ascending
central series are $\p$-ideals, there is no need to distinguish nilpotency of $(L, \p)$
from nilpotency of the underlying Lie algebra $L$.  There is likewise no need to
distinguish solubility from that of the underlying algebra.  That $(L, \p)$ has a
\p-chief series with all quotients abelian if $L$ is soluble follows from the following
lemma.

\begin{lemma}\label{lem-chab} Let $(L, \p)$ be a restricted Lie algebra whose underlying algebra $L
\ne 0$ is soluble.  Then $(L, \p)$ has a non-zero abelian \p-ideal. \end{lemma}
\begin{proof} If $\z(L) \ne 0$, then it is the required non-zero abelian \p-ideal, so
suppose $\z(L) = 0$.  Let $A$ be a minimal ideal of $L$.  For any $a \in A$, $\ad(a)^p =
0$, so $a^\scp \in \z(L) = 0$, and $A$ is a \p-ideal. \end{proof}

Note however, that the derived algebra of a restricted Lie algebra need not be a
\p-ideal.

\begin{example} \label{der} Let $L = \langle a,b,c \rangle$ with multiplication given by
$ab=b$, $ac=bc=0$.  By Jacobson's Theorem (Theorem \ref{JacThm}), $L$ has a $p$-operation
with
$a^\scp = a$, $ b^\scp = c$ and $c^\scp = 0$.  Then $L' = \langle b \rangle$ which is not a
\p-ideal.
\end{example}

\begin{example} \label{nilder}  Let $N = \langle a,b,c,d\rangle$ with multiplication given by
$ab=c$, $ac=bc=ad=bd=cd=0$.  Setting $a^\scp = b^\scp = d^\scp = 0$ and $c^\scp = d$ defines a
$p$-operation on the nilpotent algebra $N$.  Then $N' = \langle c\rangle$ which is not a
\p-ideal. \end{example}

In the following, I shall define analogues for restricted Lie algebras of concepts used
in the theory of ordinary Lie algebras.  Where these refer to a class of restricted Lie
algebras, I attach the prefix ``$p$-'' to the name of the concept.  Where it depends on
a particular $p$-operation, I attach the $p$-operation as prefix.  Thus, I shall refer
to $p$-formations and, as above, to \p-subalgebras. Where the meaning is clear from the context, I
shall often simplify notation by writing $L$ rather than $(L,\p)$. 

Much of the theory of Schunck classes and projectors relies on lemmas asserting that, under
certain circumstances, the Lie algebra $L$ splits over some abelian ideal $A$, that is,
that there exists a subalgebra $U$ such that $U+A = L$ and $U \cap A = 0$.  We shall
need lemmas giving the existence of a $\p$-subalgebra with these properties.  Note that
it is possible for the underlying Lie algebra $L$ to split over a $\p$-ideal $A$ without
there being any $\p$-subalgebra which complements $A$.

\begin{example} Let $L = \langle a,b \rangle$ with $ab=0$ and $a^\scp = 0, b^\scp = a$.  Then $A
= \langle a \rangle$ is a $\p$-ideal which is complemented in the underlying algebra but
which has no complementary $\p$-subalgebra. \end{example}

\begin{lemma}\label{auto} Let $(L, \p)$ be a restricted Lie algebra. Suppose $\z(L) =
0$.  Let $\alpha$ be an automorphism of $L$.  Then $\alpha(x)^\scp = \alpha(x^\scp)$.
\end{lemma}

\begin{proof} For $x,y \in L$,  we have
$$\ad(\alpha(x))^p \alpha(y) = \alpha(\ad(x)^p y) = \alpha(x^\scp y) = \alpha(x^\scp) \alpha(y).$$
Thus $\ad(\alpha(x))^p = \ad(\alpha(x^\scp))$.  As $\z(L) = 0$, this implies $\alpha(x)^\scp
= \alpha(x^\scp)$.  \end{proof}

\begin{lemma} \label{compl} Suppose $A$ is an abelian non-central minimal $\p$-ideal of
$(L, \p)$ and that $M$ is a subalgebra of $L$ which complements $A$.  Then $M$ is a
maximal $\p$-subalgebra of $L$. \end{lemma}

\begin{proof}   We have to prove that $M$ is a $\p$-subalgebra of $L$.  That it is
maximal then follows.  Let $x \in M$.  Then $x^\scp$ is uniquely expressible in the form
$x^\scp = x' + a$ with $x' \in M$ and $a \in A$.  For $y \in M$, $\ad(x)^p y \in M$.  But 
$$\ad(x)^p y = x^\scp y = x'y + ay,$$
so $ay \in M$.  But $ay \in A$, so $ay = 0$ for all $y \in M$.  Since $A$ is abelian, $ay
= 0$ also for all $y \in A$, so $a \in \z(L) \cap A = 0$.  Thus $x^\scp \in M$. \end{proof}

The theory of Schunck classes of soluble Lie algebras makes use of
primitive algebras.  I set out here the properties of primitive restricted Lie algebras.

\begin{definition}  A soluble restricted Lie algebra $(L, [p])$ is called
\textit{primitive} if it has a minimal $\p$-ideal $A$ with $\cent_L(A) = A$.  The minimal \p-ideal $A$ is
called the {\em socle} of $(L,\p)$ and denoted by $\soc(L,\p)$.
\end{definition}

If $A$ is an $L$-module, we put $A^L = \{a \in A \mid xa = 0
\text{ for all $x \in L$}\}$.  We use $H^n(L, A)$ to denote the ordinary cohomology of
$L$ acting on $A$.  The following lemma was proved by Voigt in \cite[Remark 2.12, p.93]{Vgt} 
using the theory of group schemes.

\begin{lemma}\label{prim}  Let $(L, \p)$ be a primitive soluble restricted Lie algebra and  let $A =
\soc(L,\p)$.  Then for all $n$, we have $H^n(L/A, A) = 0$ and there exists a subalgebra
$M$ which complements $A$, all such are conjugate under automorphisms of the form
$\alpha_a = 1 + \ad(a)$ for $a \in A$ and are maximal $\p$-subalgebras of $(L, \p)$.
\end{lemma} 

\begin{proof}  The result is trivial if $A = L$, so suppose $A \ne L$.   Let $B/A$ be a
minimal $\p$-ideal of $L/A$.  Then $A^{B/A} = 0$ and so $H^\beta(B/A, A) = 0$ for all
$\beta$.  Thus $H^\alpha(L/B, H^\beta(B/A, A)) = 0$ for all $\alpha, \beta$ and we have
$H^n(L/A, A) = 0$ for all $n$.  It follows that $L$ splits over $A$ as ordinary Lie
algebra and that all complements to $A$ in $L$ are conjugate as asserted.  That the
complements are maximal $\p$-subalgebras follows by Lemma \ref{compl}  \end{proof}

There is only one minimal soluble Lie algebra, namely the $1$-dimensional algebra with
zero multiplication.  Taking this with the zero $p$-operation gives a restricted Lie
algebra.  The following restricted Lie algebras appear in Hochschild \cite{H}, where they are called
{\em strongly abelian}.

\begin{definition} A restricted Lie algebra $(L, \p)$ is called \textit{null} if it is
abelian and $L^\scp = 0$. \end{definition}

The $1$-dimensional null algebra is not the only minimal soluble restricted Lie
algebra.  For example, for $A = \langle a \rangle$, we can set $a^\scp = a$ giving a
non-null $1$-dimensional restricted Lie algebra.  Depending on the field, there could be
other abelian algebras with no proper $\p$-subalgebras.  I shall call any of these
minimal objects an {\em atom.}

\begin{lemma}\label{chief}  Let $A/B$ be a $\p$-chief factor of the soluble restricted Lie
algebra $(L, \p)$.  Then $A/B$ is abelian.  Either $A/B$ is null, $A^\scp \subseteq B$, or
$A/B$ is central and is an atom. \end{lemma}

\begin{proof} By Lemma \ref{lem-chab}, the \p-chief factors of $(L,\p)$ are abelian.     As we can
work in $L/B$, we may suppose $B=0$.  Let $A_0 \subseteq A$ be a minimal ideal of $L$.  Suppose
$a \in A_0$ and that $a^\scp \ne 0$.  Since $\ad(a)^2 = 0$, $a^\scp \in \z(L)$.  Thus $A \cap \z(L) \ne
0$ and it follows that $A \subseteq \z(L)$.  Any $\p$-subalgebra of $A$ is a $\scp$-ideal of $L$, so 
$A$ is an atom.  So, if $A \not\subseteq \z(L)$, we have $a^\scp = 0$ for all $a \in A_0$, $A_0$ is
a null \p-ideal and $A_0 = A$.
\end{proof}

Atoms are primitive restricted Lie algebras.  Any non-abelian primitive restrictable Lie
algebra has trivial centre and so has only one $p$-operation.

\begin{lemma}\label{abpdid} Let $(L,\p)$ be a restricted Lie algebra and let $A$ be an abelian ideal     
of the underlying algebra $L$.  Then there exists a $p$-operation $\pd$ on $L$ such that $A$ is a
null \pd-ideal of $L$. \end{lemma}

\begin{proof}  Take a basis $a_1, \dots, a_r$ of $A$ and extend with elements $b_1, \dots, b_s$ to a
basis of $L$.  Since $\ad(a_i)^2=0$, by Theorem \ref{JacThm}, there exists a $p$-operation
$\pd$ on $L$ with $a_i^\scpd = 0$ and $b_j^\scpd = b_j^\scp$.  Then $A$ is a null \pd-ideal of $L$.
\end{proof}

\begin{Cor}\label{restrQ}  Let $K$ be an ideal of the soluble restrictable Lie algebra $L$.   
Then $L/K$ is restrictable. \end{Cor}

\begin{proof} If $K = 0$, then the result holds, so suppose $K \ne 0$.  Let  $A \subseteq K$ be a 
minimal ideal of $L$.  By Lemma \ref{abpdid}, $A$ is a \pd-ideal of $L$ for some $p$-operation $\pd$ on
$L$.  Thus $L/A$ is restrictable.  By induction, $L/K$ is restrictable. \end{proof}

\begin{lemma}\label{ordprim}  Let $(L, \p)$ be a soluble restricted Lie algebra and let
$L/K$ be a non-abelian primitive quotient of  $L$.  Then $K$ is a \p-ideal of
$L$. \end{lemma}

\begin{proof}  We may suppose $K \ne 0$.  Let  $A \subseteq K$ be a minimal ideal of $L$. 
By Lemma \ref{abpdid}, $A$ is a \pd-ideal of $L$ for some $p$-operation $\pd$ on $L$.  By induction,
$K/A$ is a \pd-ideal of $L/A$.  But $K \supseteq \z(L)$ since $L/K$ is non-abelian primitive.  For $k \in
K$, we have $k^\scp - k^\scpd \in \z(L)$ and $k^\scpd \in K$.  Thus $k^\scp \in K$ and $K$ is a
\p-ideal of $L$. \end{proof}

\begin{lemma}\label{ab} Suppose every \p-chief factor of $(L, \p)$ is non-null.  Then $L$
is abelian. \end{lemma}

\begin{proof}  Since every \p-chief factor is central, $L$ is nilpotent.  Let $A$ be a
minimal $\p$-ideal of $L$.  Then $A$ is central and, by induction over $\dim(L)$, $L/A$
is abelian.  For all $x \in L$, we have $\ad(x)^2 = 0$ and so $x^\scp \in \z(L)$.  But
$\langle L^\scp \rangle = L$, so $\z(L) = L$. \end{proof}

\begin{lemma}\label{nonnull} Let $A$ be a null minimal $\p$-ideal of $L$.  Suppose that
every chief factor of $L/A$ is non-null.  Then there exists a maximal $\p$-subalgebra
$M$ which complements $A$.  \end{lemma}

\begin{proof}  By Lemma \ref{ab}, $L/A$ is abelian.  Suppose $A$ is not central.  Then
$A^L = 0$, $H^n(L/A, A) = 0$ and there exists a subalgebra $M$ of $L$ which complements
$A$.  By Lemma \ref{compl}, $M$ is a maximal $\p$-subalgebra of $L$. 

Suppose $A$ is central.  By Lemma \ref{ab}, $L/A$ is abelian and it follows that
$\ad(x)^2 = 0$ for all $x \in L$ and $x^\scp \in \z(L)$.  Since $\langle (L/A)^\scp \rangle =
L/A$, we have $\z(L) + A = L$, so $L$ is abelian.  Let $M = \langle L^\scp \rangle$.  Then
$M$ is a $\p$-subalgebra of $L$ and $M + A = L$. \end{proof}

\section{The \p-Frattini subalgebra}\label{frat}

\begin{definition}  The $\p$-\textit{Frattini subalgebra} $\Psi(L, \p)$ of $(L,\p)$ is
the intersection of the maximal $\p$-subalgebras of $(L, \p)$. \end{definition}

The following is Voigt \cite[Theorem 2.88, p.247]{Vgt}.

\begin{lemma} \label{PsiId}  If $(L, \p)$ is soluble then $\Psi(L,\p)$ is a
$\p$-ideal of $L$. \end{lemma}

\begin{proof} Trivially, $\Psi$ is a $\p$-subalgebra.  We have to prove it an ideal. We
use induction over $\dim(L)$.   Let $A$ be a minimal $\p$-ideal of $L$.  By induction,
the intersection $N$ of the maximal $\p$-subalgebras which contain $A$ is an ideal.  If
every maximal $\p$-subalgebra contains some minimal \p-ideal, then the result holds,
so suppose there exists a maximal $\p$-subalgebra $M$ which does not contain any
minimal \p-ideal.  Then $M \cap A = 0$ and $M + A = L$.  If $\cent_L(A) \ne A$, then $M
\cap \cent_L(A)$ is a \p-ideal, so $\cent_L(A) = A$.  Thus $\z(L) = 0$, $A$ is a minimal
ideal of $L$ and $L$ is primitive.  The complements to $A$ are \p-subalgebras, so
$\Psi(L) = \Phi(L) = 0$.
\end{proof}

For an element $a$ of any finite-dimensional Lie algebra $L$, the Engel subalgebra
$E_L(a)$ is the Fitting null space $\ker(\ad(a)^n)$ of $\ad(a)$ for sufficiently large
$n$.  It is a subalgebra with the property that any subalgebra containing $E_L(a)$ is
self-normalising.  Also $L$ is the vector space direct sum $$\im(\ad(a)^n) \oplus
\ker(\ad(a)^n)$$ for sufficiently large $n$.

\begin{lemma}\label{engel} Let $(L, \p)$ be a restricted Lie algebra and let $a\in L$. 
Then $E_L(a)$ is a \p-subalgebra of  $(L, \p)$. \end{lemma}

\begin{proof}  If $b \in E = E_L(a)$, then $a b^\scp = - b^\scp a = - \ad(b)^p a \in E$ since
$E$ is a subalgebra.  But then $\ad(a)^n (a b^\scp) = 0$ for some $n$, that is, $b^\scp \in
E$. \end{proof}
The following result was proved by Voigt \cite[Corollary 2.92, p. 253]{Vgt} for $F$ algebraically closed.

\begin{Cor}  Let $(L, \p)$ be a restricted Lie algebra.  Suppose every maximal
\p-subalgebra of $(L, \p)$ is an ideal.  Then $(L, \p)$ is nilpotent. \end{Cor}

\begin{proof}  The \p-subalgebra $E_L(a)$ is not  contained in
any maximal \p-subalgebra.  Hence $E_L(a) = L$ for all $a \in L$ and $L$ is nilpotent by
Engel's Theorem.
\end{proof}

\begin{definition} A \p-subalgebra $S$ of $(L,\p)$ is called \textit{\p-subnormal} if
there exists a chain
$$S = S_0 \subseteq S_1 \subseteq \dots \subseteq S_n = L$$
where each $S_i$ is a \p-ideal of $S_{i+1}$. \end{definition}

A weak form of the following is given in Voigt \cite[Theorem 2.86, p. 246]{Vgt}.  (A stronger result is
given in Theorem \ref{b-n} below.)

\begin{lemma} \label{psinil}   Let $(L, \p)$ be a (not necessarily soluble) restricted Lie
algebra and let $A$ be a \p-subnormal subalgebra of $(L, \p)$.  Let $B \subseteq A \cap
\Psi(L, \p)$ be a \p-ideal of $A$.  Suppose $A/B$ is nilpotent.  Then $A$ is nilpotent.
\end{lemma}

\begin{proof}  For $a \in A$ and sufficiently large $n$, we have $\im(\ad(a)^n)
\subseteq B$ since $A$ is subnormal and $A/B$ is nilpotent.  Therefore $B + E_L(a) = L$. 
But $B \subseteq \Psi(L,\p)$, so $E_L(a) = L$.  By Engel's Theorem, $A$ is nilpotent.
\end{proof}

\begin{Cor} \label{PsiNil} Let $(L, \p)$ be a soluble restricted Lie algebra.  Then
$\Psi(L,\p)$ is nilpotent. \end{Cor}

\begin{proof} By Lemma \ref{PsiId}, $\Psi(L,\p)$ is a \p-ideal. \end{proof}

\begin{lemma} \label{psiphi}   Let $(L, \p)$ be a soluble restricted Lie algebra.  Then
$$\Psi(L,\p) \supseteq \Phi(L).$$ \end{lemma}

\begin{proof}  We use induction over $\dim(L)$.  Let $A$ be a minimal \p-ideal of $(L,
\p)$.  Then the intersection of the maximal \p-subalgebras containing $A$ contains the
intersection of the maximal subalgebras of the underlying algebra $L$ which contain
$A$.  It follows that the intersection of all maximal \p-subalgebras contains the
intersection of all maximal subalgebras unless there exists a maximal \p-subalgebra $M$
which contains no minimal \p-ideal.  But then $M$ complements the minimal \p-ideal $A$. 
Since $\cent_M(A)$ is a \p-ideal, we must have $\cent_L(A) = A$, $(L,\p)$ is primitive
and $\Psi(L,\p) = 0$.  As $A$ is also a minimal ideal of $L$, $L$ is primitive and
$\Phi(L) = 0$.
\end{proof}

Note that $\Psi(L, \p)$ can be strictly greater than $\Phi(L)$, as is the case in
Examples \ref{der}, \ref{nilder}.

\section{Schunck Classes and Projectors}\label{proj}
Following the notations of Doerk and Hawkes, I denote the class of all soluble restricted Lie algebras by
$\f{S}_p$, the class of nilpotent restricted Lie algebras by $\f{N}_p$, the clas of abelian restricted Lie
algebras by $\f{A}_p$ and the class of primitive restricted Lie algebras by $\f{P}_p$, while $\f{S},
\f{N}, \f{A}$ and $\f{P}$ denote the corresponding class of ordinary Lie algebras.  For a class
$\f{X}$, I define
\begin{equation*}\begin{split}
\quot \f{X} &= \{(L/K, \p) \mid (L, \p) \in \f{X}\}\\
\sdir \f{X} &= \{(L, \p) \mid \exists \text{ \p-ideals $K_i$ of $(L,\p)$  with $(L/K_i) \in \f{X}$
 and $\cap_i K_i = 0$}\}\\
\frat \f{X} &= \{(L,\p) \mid \exists \text{ \p-ideal $K$ of $(L,\p)$  with $K \le \Psi(L, \p)$ and $(L/K,
\p) \in \f{X}$}\}\\
\prim \f{X} &= \{(L,\p) \mid \quot(L, \p) \cap \f{P}_p \subseteq \f{X}\}.\\
\end{split}\end{equation*}
Thus $\quot\f{X}$ is the class of quotients of restricted Lie algebras in $\f{X}$, $\sdir\f{X}$ is the  
class of subdirect sums and $\frat\f{X}$ the class of Frattini extensions of algebras in $\f{X}$, while
$\prim\f{X}$ is the class of all algebras whose primitive quotients are in $\f{X}$.

\begin{definition} A non-empty class $\f{X}$ of soluble restricted Lie algebras which is $\quot$-closed,
that  is,   $\quot\f{X} =\f{X}$, is called a {\em $p$-homomorph}.  An $\sdir$-closed $p$-homomorph is
called a {\em $p$-formation.}  A non-empty class which is $\frat$-closed is called {\em saturated.}     
A non-empty class $\f{X}$ satisfying $\prim\f{X} = \f{X}$ is called a {\em $p$-Schunck class.}
\end{definition}

These definitions differ from those of Doerk and Hawkes by the inclusion of the requirement, 
convenient for the  theory of restricted Lie algebras but not for that of finite groups,  that the classes    
be non-empty.  Note also that {\em saturation} had a different meaning, explained below, in the older
terminology.  Clearly, a $p$-Schunck class is a saturated $p$-homomorph.  If $\f{X}$ is a non-empty
class, then $\prim\f{X}$ is a $p$-Schunck class, and if $\f{X}$ is a $p$-homomorph, it is the smallest
$p$-Schunck class containing $\f{X}$.

\begin{lemma}  Let $\f{X}$ be  $p$-homomorph which is not a $p$-formation.  Then there exists a
restricted Lie algebra $(L,\p)$ with \p-ideals, $K_1, K_2$ such that $(L/K_i,\p) \in \f{X}$, $i=1,2,$ but
$(L,\p) \not\in \f{X}$. \end{lemma}
\begin{proof} Doerk and Hawkes \cite[Proposition 2.5, p.272]{DH}. \end{proof} 

\begin{definition}  Let $\f{X}$ be a class of restricted Lie algebras.  A  \p-subalgebra $U$ of $(L,\p)$ is
called {\em $\f{X}$-maximal} in $(L,\p)$ if it is maximal in the set of those \p-subalgebras of $(L,\p)$
which are in $\f{X}$.
\end{definition}

\begin{definition} Let $\f{X}$ be a $p$-homomorph.  A \p-subalgebra $U$ of $(L,\p)$ is called an
{\em $\f{X}$-projector} of $(L,\p)$ if, for every \p-ideal $K$ of $L$, $U+K/K$ is $\f{X}$-maximal
in $L/K$. \end{definition}

\begin{definition} Let $\f{X}$ be a $p$-homomorph.  A \p-subalgebra $U$ of $(L,\p)$ is called an
{\em $\f{X}$-covering subalgebra} of $(L,\p)$ if, whenever $V$ is a \p-subalgebra of $L$
containing $U$ and $K$ is a \p-ideal of $V$ with $V/K \in \f{X}$, we have $U+K = V$.
\end{definition}
Thus, an $\f{X}$-covering subalgebra $U$ of $(L,\p)$ is an $\f{X}$-projector of every \p-subalgebra of
$L$ which contains $U$.  We denote the (possibly empty) set of $\f{X}$-projectors of $(L,\p)$ by
$\Proj_\f{X}(L,\p)$ and the set of $\f{X}$-covering subalgebras by $\Cov_\f{X}(L,\p)$.  

\begin{lemma}\label{Hcomp} Let $\f{X}$ be a $p$-homomorph and let $A$ be a minimal
$\p$-ideal of $(L,\p)$.  Suppose $L/A \in \f{X}$, $L \not\in \f{X}$ and that $U \in
\Proj_\f{X}(L)$.  Then $U$ complements $A$ in $L$. \end{lemma}

\begin{proof}  We have $U+A = L$.  As $U \cap A$ is a $\p$-ideal of $L$, $U \cap A =
0$. \end{proof}

Our next lemma requires a condition which is automatic for ordinary Lie algebras.

\begin{lemma} Suppose $\f{X}$ is a $p$-homomorph which contains all atoms.  Suppose $H \in
\Cov_\f{X}(L)$ and that $H$ is contained in the $\p$-subalgebra $U$ of $L$.  Then
$\sN_L(U) = U$. \end{lemma}

\begin{proof} $N = \sN_L(U)$ is a \p-subalgebra of $L$ and $U$ is a \p-ideal of
$N$.  If $N \ne U$ then there exists a \p-subalgebra $K$ of $N$ such that $K/U$ is an
atom.  Since $K/U \in \f{X}$, $H+U = K$.  But $H \subseteq  U$, so $K \subseteq U$
contrary to the choice of $K$. \end{proof}

\begin{lemma}\label{primcov} Let $\f{X}$ be a $p$-homomorph and let $(L,\p)$ be a primitive algebra
not in
$\f{X}$ but with $L/\soc(L) \in \f{X}$.  Then $\Cov_\f{X}(L) = \Proj_\f{X}(L)$ and is the set of all
complements to $\soc(L)$ in $L$. \end{lemma}

\begin{proof}  Clearly, the projectors are those maximal \p-subalgebras which do not contain
$\soc(L)$.  As $\soc(L)$ is the only minimal \p-ideal of $L$, these maximal \p-subalgebras are
easily seen to be $\f{X}$-covering subalgebras. \end{proof}

\begin{lemma} \label{conj} Let $\f{X}$ be a \p-homomorph. 
Let $A$ be an abelian \p-ideal of $(L, \p)$. Suppose $L/A \in \f{X}$ and that $U_1, U_2
\in \Cov_\f{X}(L)$.  Then there exists $a \in A$ such that $\alpha_a(U_1) = U_2$. \end{lemma}

\begin{proof}  Suppose first that $A$ is a minimal \p-ideal.  If $L \in \f{X}$, then $U_1 = U_2 = L$. 
Suppose $L \not\in \f{X}$.  Then $U_1, U_2$ complement $A$ in $L$.  The result holds if $L$ is
primitive.  If $L$ is not primitive, there exists a minimal \p-ideal $B \subseteq \cent_{U_2}(A)$
of $L$.  Since $(U_i+B)/B \in \Cov_\f{X}(L/B)$, by induction, there exists $a \in A$ such that
$$\alpha_a(U_1 + B) = U_2 + B = U_2.$$  It follows that $\alpha_a(U_1) = U_2$ since
$U_1 \simeq L/A \simeq U_2$.

Now suppose that $A_1$ is a minimal \p-ideal contained in $A$.  By induction, the result holds in
$L/A_1$, so by replacing $U_1$ by a suitable $ \alpha_a(U_1)$, we may suppose that $A_1 + U_1 =
A_1 + U_2$.  Either $A_1 + U_1\subset L$ and the result holds by induction, or $L/A_1 \in \f{X}$ and we
have the case already proved.
\end{proof}

\begin{lemma}\label{transp} Let $\f{X}$ be a $p$-homomorph.  Let $K$ be a \p-ideal of $(L,\p)$. 
Suppose $V/K \in \Proj_\f{X}(L/K)$ and $U \in \Proj_\f{X}(V)$.  Then $U \in \Proj_\f{X}(L)$.
\end{lemma}

\begin{proof} Doerk and Hawkes \cite[Proposition 3.7, p.290]{DH}. \end{proof}

\begin{lemma}\label{transc} Let $\f{X}$ be a $p$-homomorph.  Let $K$ be a \p-ideal of $(L,\p)$. 
Suppose $V/K \in \Cov_\f{X}(L/K)$ and $U \in \Cov_\f{X}(V)$.  Then $U \in \Cov_\f{X}(L)$.
\end{lemma}

\begin{proof} Doerk and Hawkes \cite[Proposition 3.7, p.290]{DH}. \end{proof}

\begin{lemma} \label{H1} Let $\f{X}$ be a \p-homomorph which contains all atoms. 
Let $A$ be a minimal \p-ideal of $(L, \p)$. Suppose $L/A \in \f{X}$, $L \not\in \f{X}$
and $\Cov_\f{X}(L) \ne \emptyset$.  Then $\Cov_\f{X}(L)$ is the set of all complements to $A$ in
$L$ and $H^1(L/A, A) = 0$.
\end{lemma}

\begin{proof}  Let $U_1 \in \Cov_\f{X}(L, \p)$.  By Lemma \ref{Hcomp},
$U_1$ complements $A$ in $L$.  If $A$ is central, then $U_1$ is a \p-ideal of $L$ and
$L/U_1 \in \f{H}$ since $\f{H}$ contains all atoms.  Therefore $A$ is not
central in $L$. Thus $A$ is null and is a minimal ideal of $L$.  Let $U_2$ be
another complement to $A$ in $L$.  Then $(U_2, \p) \simeq (L/A, \p) \simeq (U_1, \p)$,
so $U_2 \in \f{X}$.  We have to prove that $U_2 \in \Cov_\f{H}(L)$. 

Suppose $\z(L) = 0$.  There is an automorphism $\alpha$ of the underlying algebra $L$
which maps $U_1$ onto $U_2$.  By Lemma \ref{auto}, $\alpha$ is an automorphism of $(L,
\p)$, and it follows that $U_2 \in \Cov_\f{X}(L)$.  Now suppose $\z(L) \ne 0$. 
Let $Z \subseteq \z(L)$ be a minimal \p-ideal of $L$.  Since $U_i$ is not an ideal of
$L$, we have $U_i \supset Z$.  By induction, $U_2/Z \in \Cov_\f{X}(L/Z)$.  As $U_2 \in
\f{H}$, by Lemma \ref{transc}, $U_2 \in \Cov_\f{X}(L)$.

By Lemma \ref{conj}, it now follows for every complement $U_2$ to $A$ in $L$, that there
exists $a \in A$ such that $U_2 = \alpha_a(U_1)$.  This is equivalent to $H^1(L/A, A) =
0$.
\end{proof}

\begin{definition} A $p$-homomorph $\f{X}$ is called {\em projective} if, for every soluble restricted
Lie algebra $(L,\p)$, we have $\Proj_\f{X}(L) \ne \emptyset$.  It is called a {\em Gasch\"utz} class
if, for every soluble restricted Lie algebra $(L,\p)$, $\Cov_\f{X}(L) \ne \emptyset$. \end{definition} 

In the older terminology, what are here called $\f{X}$-covering subalgebras were called
$\f{X}$-projectors, and what are here called Gasch\"utz classes were called saturated homomorphs.

\begin{theorem} \newcounter{bean} Let $\f{X}$ be a $p$-homomorph.  Then the following are
equivalent:
\begin{list}{}{\usecounter{bean}}
\item[\rm (a)] $\f{X}$ is a Gasch\"utz class.
\item[\rm (b)] $\f{X}$ is a projective class.
\item[\rm (c\,)] $\f{X}$ is a $p$-Schunck class. \end{list}\end{theorem}

\begin{proof} If $\f{X}$ is a Gasch\"utz class, it is clearly a projective class.  Suppose $\f{X}$ is
projective.  To show that it is a $p$-Schunck class, we have to how that, if $(L,\p) \in \prim\f{X}$, the
$(L,\p) \in \f{X}$.  Let $(L,\p)$ be a minimal counterexample and let $A$ be a minimal \p-ideal of $L$. 
Then $L/A \in\f{X}$.  Since $\f{X}$ is projective, there exists an $\f{X}$-projector $U$ of $L$.  If
there exists a minimal \p-ideal $B$ of $L$ contained in $U$, then $L/B \in \f{X}$ contrary to $U$ being
an $\f{X}$-projector.  Therefore $\cent_L(A) = A$ and $L$ is primitive.  But by assumption, all
primitive quotients of $L$ are in $\f{X}$, so $L \in \f{X}$ contrary to assumption.

Now suppose that $\f{X}$ is a $p$-Schunck class.  We use induction over the dimension to show that
every soluble restricted Lie algebra $(L,\p)$ has an $\f{X}$-covering subalgebra.  We may suppose $L
\not\in
\f{X}$.  Thus there exists a \p-ideal with $P = L/K$ primitive, not in $\f{X}$ but with $P/\soc(P) \in
\f{X}$.  If $K \ne 0$, then by induction, there exists an $\f{X}$-covering subalgebra $U/K$ of $L/K$. 
As $U \ne L$, by induction, we have there exists an $\f{X}$-covering subalgebra $V$ of $U$.  By
Lemma \ref{transc}, $V \in \Cov_\f{X}(L)$.  If $K=0$, then $L$ is primitive and  the complements to
$\soc(L)$ are $\f{X}$-covering subalgebras of $L$ by Lemma  \ref{primcov}. \end{proof}

\begin{theorem} \label{anyop} Let $\f{H}$ be a $p$-Schunck class which contains all
atoms.  Let $(L, \p) \in \f{H}$ and let \pd\ be another $p$-operation on $L$.  Then $(L,
\pd) \in \f{H}$. \end{theorem}

\begin{proof} Put $Z = \z(L)$.   Let $K$ be a \pd-ideal of $L$ with $L/K$ primitive.  
If $K \supseteq Z$, then $(L/K,\pd) = (L/K, \p) \in \f{H}$.  If $K
\not\supseteq Z$, then $(L/K, \pd)$ is an atom and so in $\f{H}$.  Thus $(L,
\pd) \in  \f{H}$. \end{proof}

\begin{lemma}\label{closure} Let $\f{H}$ be a $p$-homomorph.  Suppose $U$ is an $\f{H}$-covering
subalgebra of $(L,\p)$.    Then $U$ is a  $\prim\f{H}$-covering subalgebra of $(L,\p)$.
 \end{lemma}

\begin{proof}  Let $(L, \p)$ be a minimal counterexample.  Let $A$ be a minimal \p-ideal of $L$. 
Then $U+A/A$ is a $\prim\f{H}$-covering subalgebra of $(L/A)$.  If $U+A \subset L$, then $U$ is a
$\prim\f{H}$-covering subalgebra of $U+A$ and so also of $L$ by Lemma \ref{transc}.  Hence $U+A
= L$ and $U$ complements $A$ in $L$.  As this holds for every minimal \p-ideal, $U$ contains no
non-trivial \p-ideal of $L$ and  it follows that $\cent_L(A) = A$.  As $L$  is primitive and
not in $\f{H}$, $L \not\in \prim{\f{H}}$.  But $L/A \in \prim{\f{H}}$ and the
complements to $A$ in $L$ are $\prim{\f{H}}$-covering subalgebras. \end{proof}

Our next lemma is a slightly modified version of Doerk and Hawkes \cite[Lemma 3.14, p.295]{DH}.
\begin{lemma} \label{nilsup} Let $\f{X}$ be a $p$-Schunck class.  Let $N$ be a nilpotent \p-ideal
of $(L,\p)$ and  let $U$ be an $\f{X}$-maximal \p-subalgebra of $L$ such that $L = U + N$.  Then $U
\in \Cov_\f{X}(L)$. \end{lemma}

\begin{proof}  We use induction on the dimension of $L$.  If $L \in \f{X}$, the result holds, so suppose
$L \not\in \f{X}$.  Then there exists a \p-ideal $K$ of $L$ such that $P=L/K$ is primitive, not in
$\f{X}$ but with $P/\soc(P) \in \f{X}$.  Then $N \not\subseteq K$ since $L/N \simeq U/(U\cap N) \in
\f{X}$.  Hence $N+K/K$ is a non-zero nilpotent \p-ideal of $P$ and so $N+K/K = \soc(P)$.  But $U+K/K$
is a maximal \p-subalgebra of $P$ complementing $N+K/K$, so $U+K/K \in \Cov_\f{X}(L/K)$ by
Lemma \ref{primcov}.  Now $(U+K) \cap N$ is a nilpotent \p-ideal of $U+K$. Also $U+((U+K) \cap N) =
U+K$ by the modular law for subspaces.  As $U+K \ne L$, by induction, we have that $U \in
\Cov_\f{X}(U+K)$.  By Lemma \ref{transc}, $U \in \Cov_\f{X}(L)$. \end{proof}

\begin{theorem}\label{projcov} Let $\f{X}$ be a $p$-Schunck class.  Let $U \in \Proj_\f{X}(L,\p)$. 
Then $U \in \Cov_\f{X}(L,\p)$. \end{theorem}

\begin{proof} Let $A$ be a minimal \p-ideal of $L$ and put $V = U+A$.  Then $V/A \in
\Proj_\f{X}(L/A)$, so by induction, $V/A \in \Cov_\f{X}(L/A)$.  But $U$ is $\f{X}$-maximal in $V$,
$A$ is a nilpotent \p-ideal of $V$ and $U+A=V$.  By Lemma \ref{nilsup}, $U \in \Cov_\f{X}(V)$.  By
Lemma \ref{transc}, $U \in \Cov_\f{X}(L)$. \end{proof}

\begin{lemma}\label{Fcomp} Let $\f{F}$ be a $p$-formation and let $A$ be a minimal
$\p$-ideal of $L$.  Suppose $L/A \in \f{F}$, $L \not\in \f{F}$ and that $H$ complements
$A$ in $L$.  Then $H \in \Cov_\f{F}(L)$. \end{lemma}

\begin{proof} $H \in \f{F}$ and is a maximal $\p$-subalgebra of $L$.  Thus we have only
to prove that $K$ a $\p$-ideal of $L$, $L/K \in \f{F}$ implies $H + K = L$.  Since
$\f{F}$ is a $p$-formation, $L/K \cap A \in \f{F}$.  Since $A$ is minimal and $L \not\in
\f{F}$, this implies $K \supseteq A$ and $H+K \supseteq H + A = L$. \end{proof}

Taking $\f{F}$ to be the $p$-homomorph of null algebras of dimension at most $1$ and $L$
the $2$-dimensional null algebra shows that the assumption that $\f{F}$ is a
$p$-formation, not merely a $p$-homomorph, cannot be omitted from Lemma \ref{Fcomp}.

\begin{lemma}\label{schunck} Let $\f{H}$ be a $p$-homomorph.  A necessary and
sufficient condition for $\f{H}$ to be a $p$-Schunck class is that $L \not\in \f{H}$, $A$ a minimal
$\p$-ideal of $L$ and $L/A \in \f{H}$ implies $\f{H}(L) \ne \emptyset$. \end{lemma}

\begin{proof}  The condition is trivially necessary.  Suppose $\f{H}$ satisfies the
condition.  We use induction over $\dim(L)$ to prove for all $(L,\p)$ that $\Cov_\f{H}(L) \ne
\emptyset$.  Let $A$ be a minimal $\p$-ideal of $L$.  Then there exists a
$\p$-subalgebra $U \supseteq A$ such that $U/A \in \Cov_\f{H}(L/A)$.

If $U \subset L$, then by induction, there exists $H \in \Cov_\f{H}(U)$ and then $H \in \Cov_\f{H}(L)$
by Lemma \ref{transc}.  If $U = L$, then $L/A \in \f{H}$.  Either $L \in \f{H}$ or $L
\not\in \f{H}$ and $\f{H}(L) \ne \emptyset$ by hypothesis. \end{proof}

\begin{Cor}\label{split}  Let $\f{F}$ be a $p$-formation.  Then $\f{F}$ is a $p$-Schunck class if
and only if $\f{F}$ is saturated. \end{Cor}

\begin{proof} Suppose $\f{F}$ is a $p$-Schunck class.  Suppose $L/\Psi(L) \in \f{F}$.  We have to
prove that $L \in \f{F}$.  We may suppose $\Psi(L) \ne )$.  Let $A \subseteq \Psi(L)$ be a minimal
\p-ideal of $L$.  By induction, $L/A \in \f{F}$.  Let $U \in \Proj_\f{F}(L)$.  If $L \not\in \f{F}$, then by
Lemma \ref{Hcomp}, $U$ complements $A$ in $L$ contrary to $A \subseteq \Psi(L)$.  Therefore    $L
\in \f{F}$.

Suppose $\f{F}$ is saturated.  Suppose $(L,\p) \not\in \f{F}$ and that $A$ is a minimal \p-ideal of $l$
with $L/A \in \f{F}$.  Since $\f{F}$ is saturated, $A \not\subseteq \Psi(L)$, so there exists a maximal
\p-subalgebra $U$ complementing $A$ in $L$.  Suppose $K$ is a \p-ideal of $L$ and that $U+K/K$ is
not $\f{F}$-maximal in $L/K$.  Then we must have $K \subseteq U$ an $L/K \in \f{F}$.  But then
$L/(K \cap A) \in \f{F}$ since $\f{F}$ is a $p$-formation, that is, $L \in \f{F}$.
\end{proof}

\begin{lemma}\label{satH} Let $\f{H}$ be a $p$-Schunck class.  Then
\begin{enumerate}
\item $L,M \in \f{H}$ implies $L \oplus M \in \f{H}$, and 
\item $L/ \Psi(L) \in \f{H}$ implies $L \in \f{H}$. \end{enumerate}
\end{lemma}

\begin{proof} As for Barnes and Gastineau-Hills \cite[Lemma 3.5]{BGH}. \end{proof}

\begin{Cor}\label{central}   Let $\f{H}$ be a $p$-Schunck class  which contains all
atoms.  Suppose $A$ is a central \p-ideal of $L$ and $L/A \in \f{H}$.  Then $L \in
\f{H}$. \end{Cor}

\begin{proof}  We may suppose $A$ minimal.  If $L \not\in \f{H}$, then $L$ splits over
$A$ and we have $L \simeq (L/A) \oplus A$. \end{proof}

If a $p$-Schunck class is non-zero, it must contain some atom.  However, it need
not contain every atom.

\begin{example}  Let $\f{A}$ be the class of all those soluble restricted Lie algebras all of
whose $\p$-chief factors are non-null.  Then $\f{A}$ is clearly a formation.  By Lemma
\ref{nonnull}, $\f{A}$ is saturated.  By Lemma \ref{ab}, every algebra in $\f{A}$ is
abelian. \end{example}

The following result, which generalises Lemma \ref{psinil}, is the analogue of Barnes
and Newell \cite[Theorem 4.3]{BN}.  Omitting reference to the $p$-operation gives an
improved proof of that result.

\begin{theorem} \label{b-n} Let $\f{H}$ be a $p$-Schunck class which contains all
atoms.  Let $A$ be \p-subnormal in the (not necessarily soluble) restricted Lie algebra
$(L, \p)$ and let $B$ be a \p-ideal of $A$.  Suppose $B \subseteq \Psi(L)$ and $A/B \in
\f{H}$.  Then $A \in \f{H}$. \end{theorem}

\begin{proof}  Let $(L, \p)$ be a minimal counterexample.  If $C$ is any non-zero
\p-ideal of $L$, then $A+C/C \in \f{H}$.  We prove first that there exists a null minimal
\p-ideal $C$ of $L$ with $C \subseteq A$.  If $\z(L) \ne 0$, then $A/A \cap \z(L) \in
\f{H}$ and it follows by Corollary \ref{central}, that $A \in \f{H}$.  So $\z(L) = 0$. Since $A \not\in
\f{H}$, $A$ is not nilpotent.  By Schenkman's Theorem (see \cite{Sch}), the nilpotent residual
$A_{\f{N}}$ of $A$ is an ideal of $L$.  Let $C \subseteq A_\f{N}$ be a minimal ideal of
$L$.  By Lemma \ref{psinil}, $B$ is nilpotent, so $A$ is soluble.  Hence $C$ is
abelian.  As $\ad(c)^2 = 0$ for all $c \in C$, $c^\scp \in \z(L) = 0$. Thus $C$ is a null
\p-ideal of $L$.  

We have $A/C \in \f{H}$ but $A \not\in \f{H}$, so there exists a
\p-ideal $K$ of $A$ such that $P = A/K$ is primitive and not in $\f{H}$.  As $A/C \in
\f{H}$, $C \not\subseteq K$.  Thus $C+K/K$ is a non-zero  abelian \p-ideal of $A/K$.  As
$A/K$ is primitive, $C+K/K$ is its minimal \p-ideal. Working with the underlying
algebras, we have $M = C \cap K$ is a maximal $A$-submodule of $C$, $A$ acts
non-trivially on $C/M$, $A/M$ splits over $C/M$ and, by Lemma \ref{H1}, $H^1(A/C, C/M) =
0$.  By Barnes and Newell \cite[Lemma 4.2]{BN}, $L$ splits over $C$.   Let $U$
complement $C$ in $L$.  By Lemma \ref{compl}, $U$ is a maximal \p-subalgebra of $L$. 
Since $B \subseteq
\Psi(L)$, $U \supseteq B$.  Put $D = M+B$.
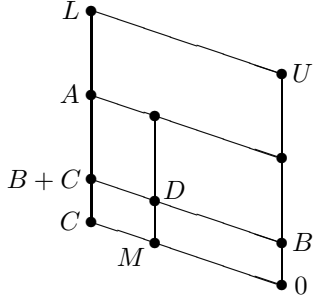
\begin{center}
\setlength{\unitlength}{0.8em}
\begin{picture}(16,16)(-7,0)
\put(-3,15){\circle*{.5}}
\put(-4.4,14.6){$L$}
\put(-3,15){\line(3,-1){9}}
\put(-3,15){\line(0,-1){10}}
\put(-3,11){\circle*{.5}}
\put(-4.5,10.6){$A$}
\put(-3,7){\circle*{.5}}
\put(-7.0,6.6){$B+C$}
\put(-3,5){\circle*{.5}}
\put(-4.5,4.6){$C$}
\put(6,12){\circle*{.5}}
\put(6.5,11.6){$U$}
\put(-3,11){\line(3,-1){9}}
\put(6,8){\circle*{.5}}
\put(-3,7){\line(3,-1){9}}
\put(6,4){\circle*{.5}}
\put(6.5,3.6){$B$}
\put(-3,5){\line(3,-1){9}}
\put(6,12){\line(0,-1){10}}
\put(6,2){\circle*{.5}}
\put(6.6,1.6){$0$}
\put(0,10){\circle*{.5}}
\put(0,6){\circle*{.5}}
\put(0,4){\circle*{.5}}
\put(0,10){\line(0,-1){6}}
\put(0.4,6.1){$D$}
\put(-1.8,2.9){$M$}
\end{picture}
\end{center} 

Then $D$ is a \p-ideal of $A$, $B+C/D$ is a complemented \p-chief factor of $A$
isomorphic to $C/M$. Thus $A/D$ has a primitive quotient isomorphic to $P$.  But this is
a primitive quotient of $A/B \in \f{H}$, contrary to $P \not\in \f{H}$.
\end{proof}

\begin{definition} Let $(L, \p)$ be a restricted Lie algebra and let $V$ be an
$L$-module giving the representation $\rho$.  $V$ is called a \textit{\p-module} for
$L$ or an \textit{$(L, \p)$-module} and $\rho$ a \textit{\p-representation} of $L$ if
$\rho(x^\scp) = \rho(x)^p$ for all $x \in L$.
\end{definition}

Clearly, an $L$-submodule of an $(L,\p)$-module
is an $(L, \p)$-submodule.  Thus, an irreducible $(L,\p)$-module is also irreducible as
$L$-module.

If $A$ is an $(L,\p)$-module, we can form the split extension of $A$ by $(L, \p)$.  This
is the ordinary split extension $X$ of $A$ by $L$ with $p$-operation coinciding with the
given operation on $L$ and null on $A$.  By Theorem \ref{JacThm}, there
is one and only one such $p$-operation.  Sometimes we are given a
$p$-operation on the module $A$ as well as on $L$, for example if $A$ is a \p-ideal of
$L$, and require the
$p$-operation on $X$ to agree with these.  For $a \in A$, $\ad_X(a)^2 = 0$, so for this
to be possible we must have $a^\scp \in A^L$.  If
this condition is satisfied, then we have $\ad_X(a)^p = \ad(a^\scp)$.  It then follows by
Theorem \ref{JacThm} that the given $p$-operation \p\ on $L$ and $A$ has a unique
extension to a $p$-operation on $X$.

\begin{definition}  A $p$-homomorph $\f{H}$ is said to be \textit{split} if $L \in \f{H}$
and $A$ an abelian \p-ideal of $L$ imply that the split extension of $A$ by $L/A$ is
also in $\f{H}$. \end{definition}

\begin{lemma}\label{Fissplit} Let $\f{F}$ be a $p$-formation.  Then $\f{F}$ is split.
\end{lemma}
\begin{proof} As for Barnes and Gastineau-Hills \cite[Lemma 1.16]{BGH}. \end{proof}

\begin{lemma}\label{SplitF} Let $\f{H}$ be a $p$-Schunck class which is split.  Then
$\f{H}$ is a $p$-formation. \end{lemma}

\begin{proof} As for Barnes and Gastineau-Hills \cite[Theorem 2.8]{BGH}. \end{proof}

\section{$\f{F}$-hypercentral modules}\label{hyperc}

\begin{definition} Let $\f{F}$ be a saturated $p$-formation and let $V$ be an irreducible
$(L,\p)$-module. $V$ is called {\em $\f{F}$-central} if the split extension of $V$ by
$L/ \cent_L(V) \in \f{F}$ and {\em $\f{F}$-eccentric} otherwise. An $(L,\p)$-module
$V$ is called {\em $\f{F}$-hypercentral} if every composition factor
of $V$ is $\f{F}$-central.  An $(L,\p)$-module $V$ is called
{\em $\f{F}$-hypereccentric} if every composition factor of $V$ is
$\f{F}$-eccentric. \end{definition}

If $A, B$ are \p-ideals of $(L, \p)$ and $A/B$ is a null \p-chief factor of $L$, then
$A/B$ is an irreducible $(L, \p)$-module and it may be classified as $\f{F}$-central or
$\f{F}$-eccentric as above.  I extend the definitions to apply to any \p-chief factor.

\begin{definition}   Let $\f{F}$ be a saturated $p$-formation which contains the null
atom.  A \p-chief factor $A/B$ of $(L, \p)$ is called {\em $\f{F}$-central} if the 
split extension of $A/B$ by $L/ \cent_L(A/B) \in \f{F}$ and {\em $\f{F}$-eccentric}
otherwise. \end{definition}

\begin{lemma}\label{allFc}  Let $\f{F}$ be a saturated $p$-formation which contains all
atoms.  Then $(L, \p) \in \f{F}$ if and only if every \p-chief factor of $L$ is
$\f{F}$-central. \end{lemma}

\begin{proof} Suppose $(L, \p) \in \f{F}$.  Let $A/B$ be a \p-chief factor.  By Lemma
\ref{Fissplit}, the split extension $X$ of $A/B$ by $L/A$ is in $\f{F}$.  Let $U \simeq
L/A$ be a complement to $A/B$ in $X$ and let $C = \cent_U(A/B)$.  Then $C$ is a \p-ideal
of $X$ and $X/C$ is the split extension of $A/B$ by $L/\cent_L(A/B)$ and is in $\f{F}$.

Suppose conversely, that every \p-chief factor of $(L, \p)$ is $\f{F}$-central.  Let $A$
be a minimal \p-ideal of $L$.  By induction, we may suppose that $(L/A, \p) \in \f{F}$. 
If $(L/A, \p) \not\in \f{F}$, then by Corollary \ref{split}, there exists a
\p-subalgebra $M$ of $L$ which complements $A$.  Let $C = \cent_M(A)$.  Then $L/C$ is
the split extension of $A$ by $L/\cent_L(A)$ and is in $\f{F}$.  Thus $L/C$ and $L/A$
are in $\f{F}$, so $L/(C \cap A) \in \f{F}$.  But $C \cap A = 0$.
\end{proof}

If $V, W$ are $(L,\p)$-modules, then $V \otimes_F W$ and $\Hom_F(V, W)$ are also
$(L, \p)$-modules.  The following results are proved exactly as the corresponding results
for ordinary Lie algebras.

\begin{theorem}\label{tensHom} Let $\f{F}$ be a saturated $p$-formation and let $V, W$
be $\f{F}$-hypercentral $(L, \p)$-modules.  Then $V \otimes_F W$ and $\Hom_F(V, W)$ are 
also $\f{F}$-hypercentral. \end{theorem}

\begin{proof}  The argument for Theorem 2.1 of \cite{HyperC} applies. \end{proof}

If $S$ is a \p-subalgebra of $(L, \p)$ and $V$ is an $(L, \p)$-module, then it is also a
$(S, \p)$-module.  Let $\f{F}$ be a saturated $p$-formation.  We say that $V$ is 
$S\f{F}$-hypercentral if it is $\f{F}$-hypercentral as $S$-module and
$S\f{F}$-hyperexcentric if it is $\f{F}$-hyperexcentric as $S$-module.
 
\begin{theorem} \label{compon}  Let $(L, \p)$ be a (not necessarily soluble) restricted Lie
algebra. Let $\f{F}$ be a saturated $p$-formation.  Suppose $S$ is \p-subnormal in $(L,
\p)$ and that $S \in \f{F}$.   Let $V$ be a finite-dimensional $(L, \p)$-module.  Then
$V$ is the $L$-module direct sum $V = V_0 \oplus V_1$ where $V_0$ is
$S\f{F}$-hypercentral and $V_1$ is $S\f{F}$-hyperexcentric.
\end{theorem}

\begin{proof}  The argument for Lemma 1.1 of \cite{extras} applies. \end{proof}

\begin{theorem} \label{thm-prep}  Let $(L, \p)$ be a (not necessarily soluble) restricted Lie
algebra and suppose that $z^\scp = 0$ for all $z \in \z(L)$.  Let $\f{F}$ be a
saturated $p$-formation and suppose $S$ is \p-subnormal in $L$, $S \ne 0$ and that
$S \in \f{F}$.  Let $V$ be an irreducible $(L, \p)$-module.  Then $V$ is
$S\f{F}$-hypercentral.
\end{theorem}

\begin{proof} The argument for Theorem 6.4 of \cite{extras} applies. \end{proof}

\section{Constructions and examples}\label{constr}
Starting with any non-empty class $\f{X}$, we can construct a $p$-Schunck  class by forming the
class $\prim(\quot \f{X})$ of all souble restricted Lie algebras whose primitive quotients are quotients
of algebras in $\f{X}$.  In this section, I give some constructions for $p$-Schunck classes with the extra
property of being formations. 

\begin{definition} Let $\f{K}$ be a saturated $p$-formation which contains all atoms and
let $\f{F}$ be a $p$-formation.  The $p$-formation \textit{residually defined}
by $\f{K}$ and $\f{F}$ is the class
$$\f{K} \cdot \f{F} = \{(L, \p) \mid L_{\f{F}} \in \f{K} \}.$$ \end{definition}

This is {\em not} the product $\f{K}\f{F}$ as defined in Doerk and Hawkes \cite[Definition 1.3,
p. 263]{DH}.  An algebra in $\f{K} \cdot \f{F}$ is an extension of an algebra in $\f{K}$ by an algebra
in $\f{F}$.  If $\f{K}$ is \p-ideal closed as is, for example, the class $p\f{N}$ of
nilpotent restricted algebras, then every restricted algebra in $\f{K}\f{F}$, that is, every extension
of an algebra in $\f{K}$ by one in $\f{F}$, is in $\f{K} \cdot \f{F}$.

\begin{theorem} \label{resdef}  Let $\f{K}$ be a saturated $p$-formation which contains all
atoms and let $\f{F}$ be a $p$-formation.  Then $\f{K} \cdot \f{F}$ is a
saturated formation. \end{theorem}

\begin{proof}  Suppose $(L, \p) \in \f{K} \cdot \fF$ and let $K$ be a \p-ideal of $L$. 
Then
$$(L/K)_{\f{F}} = (L_{\f{F}} + K)/K \simeq L_{\f{F}}/(L_{\f{F}} \cap K) \in \f{K},$$
so $(L/K, \p) \in  \f{K} \cdot \f{F}$.  Now suppose $A_1, A_2$ are \p-ideals of $(L, \p)$
and that $(L/A_i, \p) \in \f{K} \cdot \f{F}$.  We have to show that $L/(A_1 \cap A_2) \in
\f{K} \cdot \f{F}$.  We may suppose $A_1 \cap A_2 = 0$.  We have $L_{\f{F}} / (L_{\f{F}}
\cap A_i) \simeq (L_{\f{F}} + A_i)/A_i \in \f{K}$.  Thus $L_{\f{F}} \in \f{K}$ and so
$\f{K} \cdot \f{F}$ is a formation.

Now suppose $A$ is a minimal \p-ideal of $L$, $A \subseteq \Psi(L)$ and that
$L/A \in \f{K} \cdot \f{F}$.  We have to prove that $L_\f{F} \in \f{K}$.  If $A
\not\subseteq L_{\f{F}}$, then $L_{\f{F}} \simeq (L_{\f{F}}+ A)/A = (L/A)_{\f{F}} \in
\f{K}$.
   
Suppose $A \subseteq L_{\f{F}}$. Then $L_{\f{F}} / A \in \f{K}$.  If $A \subseteq
\z(L)$, then $L_\f{F} \in \f{K}$, so we may suppose that $A$ is null and therefore is a
minimal ideal of the underlying algebra $L$.  Suppose $L_\f{F} \not\in \f{K}$.     Since $A$ is
irreducible as $L$-module and $L_\f{F}$ is an ideal of $L$, all composition factors of $A$ as
$L_\f{F}$-module are isomorphic.  Let
$$A = A_0 \supset A_1 \supset \dots \supset A_n = 0$$
be a composition series of $A$ as $L_\f{F}$-module.  Since $L_\f{F} \not\in \f{K}$,
there exists $r$ such that $L_\f{F}/A_r \in \f{K}$ but $L_\f{F}/A_{r+1} \not\in \f{K}$. 
Therefore $L_\f{F}/A_{r+1}$ splits over $A_r/A_{r+1}$ and it follows that the split
extension of $A_r/A_{r+1}$ by $L_\f{F}/A$ is not in $\f{K}$.  But $\f{K}$ is a saturated
formation, so this implies that the split extension is the only extension, so
$H^2(L_\f{F}, A_r/A_{r+1}) = 0$.  We also have $H^0(L_\f{F}, A_r/A_{r+1}) = H^1(L_\f{F},
A_r/A_{r+1}) = 0$.  Since all the $A_i/A_{i+1}$ are isomorphic, we have $H^n(L_\f{F},
A_i/A_{i+1}) = 0$ for all $i$ and $n \le 2$.  It follows that $H^n(L_\f{F}, A)
= 0$ for $n \le 2$.  By the Hochschild-Serre spectral sequence, we have $H^n(L, A) = 0$
for $n \le 2$.  Thus $L$ splits over $A$.  By Lemma \ref{compl}, the complements are maximal
\p-subalgebras contrary to $A \subseteq \Psi(L)$.
\end{proof}

In particular, the class \pCs\ of completely soluble restricted Lie algebras, that is the
algebras with nilpotent derived algebras, is a saturated $p$-formation since the
\p-closure $(L')_\scp$ of  $L'$ is the \p-abelian residual and, by
Strade and Farnsteiner \cite[Proposition 1.3(3), p. 66]{SF}, $(L')_\scp$ is nilpotent if
$L'$ is nilpotent.

Likewise, by induction over $k$, it is easily seen that $p\f{N}^k$, the $p$-formation of
algebras of nilpotent length at most $k$, is saturated.

Since every \p-chief factor of a soluble restricted Lie algebra $(L, \p)$ is either a
chief factor of the underlying algebra $L$ or is central, the intersection of the
centralisers of the \p-chief factors of $(L, \p)$ is the intersection of the
centralisers of the chief factors of $L$, that is, the nil radical $N(L)$.  Let $\f{F}$
be a $p$-formation.  The $p$-formation $\pLoc(\f{F})$ $p$-locally defined by $\f{F}$ is
the class of all $(L, \p)$ such that, for every \p-chief factor $C$ of $(L, \p)$, $(L/
\cent_L(C), \p) \in \f{F}$, that is, all $(L, \p)$ with $(L/N(L), \p) \in \f{F}$ or
equivalently, $L_\f{F} \in p\f{N}$.  Thus $\pLoc(\f{F}) = p\f{N} \cdot \f{F}$ and is
saturated.

\begin{definition}  The algebra $(L, \p)$ is called \textit{$p$-supersoluble} if every
\p-chief factor of $(L, \p)$ is an atom. \end{definition}

This is equivalent to the definition given in Voigt \cite[Definition 2.1, p.77]{Vgt}.  Voigt proves
(\cite[Satz 2.17, p.117]{Vgt}), for a soluble restricted Lie algebra $(L, \p)$ over an algebraically
closed field $F$ , that $(L, \p)$ is supersoluble if and only if every maximal \p-subalgebra has
codimension $1$.  The assumption that $F$ is algebraically closed, made necessary by the existence of
atoms of dimension greater than $1$,  can be dropped if we replace ``dimension'' by ``\p-dimension'',
defined for a soluble restricted Lie algebra as the length of a \p-composition series.   Voigt's Theorem
2.18, that a soluble restricted Lie algebra is supersoluble if and only if all maximal chains of
\p-subalgebras have the same length, also  follows without the assumption of algebraic closure.

As every \p-chief factor is either a chief factor of $L$ or central, $(L, \p)$ is
$p$-supersoluble if and only if $L$ is supersoluble.  We denote the class of
$p$-supersoluble restricted Lie algebras by $p\f{U}$.

\begin{lemma} \label{supsol} $p\f{U}$ is a saturated $p$-formation. \end{lemma}

\begin{proof}  $p\f{U}$ is obviously a $p$-formation.  Suppose $A$ is a minimal \p-ideal
of $(L, \p)$, $(L/A, \p) \in p\f{U}$ and $A$ is not an atom.  Then $A$ is not central,
so is a minimal ideal of $L$.  As $L/A$ is supersoluble and $A$ has dimension greater
than $1$, $L$ splits over $A$ and the result follows. \end{proof}

\begin{definition} Let $\Lambda$ be an $F$-subspace of the algebraic closure $\bar{F}$. 
We say that $\Lambda$ is \textit{$p$-normal} if $\lambda \in \Lambda$ implies that every
conjugate of $\lambda$ is in $\Lambda$ and that $\lambda^p \in \Lambda$. \end{definition}

If $F'$ is a normal extension field of $F$, then $F'$ is a $p$-normal subspace of
$\bar{F}$.  If $p > 2$ and $k \in F$ has no square root in $F$, then the space
$\langle \sqrt{k}\ \rangle$ is another example.  If $\lambda^p = k \in F$, $\lambda \not
\in F$, then $\langle \lambda \rangle$ is a normal subspace of $\bar{F}$ which is not
$p$-normal. $F$ being perfect does not ensure that a normal space will be $p$-normal. 
The following example is based on an idea provided by G.~E.~Wall.

\begin{example} \label{notpn} Let $F$ be the field of $p^n$ elements.  Let $q > p$ be a
prime dividing $p^n - 1$.  Then $F$ contains a primitive $q$-th root $\xi$ of unity and
also has an element $c$ which has no $q$-th root in $F$.  Let $u \in \bar{F}$ be a root
of the polynomial $f(t) = t^q - c$ and let $\Lambda$ be the space spanned by the
conjugates of $u$.  Then the $u \xi^i$ for $i = 0, 1, \dots, q-1$ are the roots of
$f(t)$.  If $m(t) = t^k + a_1 t^{k-1} + \dots + a_k$ is the minimum polynomial of $u$,
then $m_i(t) = t^k + a_1\xi^i t^{k-1} + \dots + a_k (\xi^i)^k$ is the minimum polynomial
of $u \xi^i$.  Each of these divides $f(t)$.  As $f(t)$ is a product of irreducible
polynomials of the same degree $k$, the degree of $f(t)$ is divisible by $k$.  As the
degree $q$ is prime, it follows that $f(t)$ is irreducible and is the minimal polynomial
of $u$.  We therefore have $\Lambda = \langle u, \xi u, \dots \xi^{q-1} u \rangle =
\langle u \rangle$ since $\xi \in F$.  If $u^p \in \Lambda$, then $u^p = bu$ for some $b
\in F$, and we have $u^{p-1} - b = 0$.  But $t^{p-1} - b$ is a polynomial of degree less
than that of the minimum polynomial of $u$.  Hence $u^p \not\in \Lambda$. \end{example}

Such $p, n, q$ do exist, for example, $p=2$, $n=2$, $q=3$ and $p=3$, $n=3$, $q=13$. 
Indeed, for any $p$, there exist such $n$ and $q$.

\begin{lemma}  Let $p$ be prime.  Then there exists $n$ and a prime $q > p$ which divides
$p^n-1$. \end{lemma}

\begin{proof}  I show that the set $S$ of primes $q$ which divide $1 + p + \dots +
p^{k-1}$ for some $k$ is infinite and so has a member greater than $p$.  For each $q \in
S$, there is a least $k$ for which $q$ divides $1 + p + \dots + p^{k-1}$. 
Consider $n > k$.  Then
$$1+ \dots + p^{n-1} = (1+ \dots + p^{k-1}) + p^k(1+ \dots + p^{k-1}) + \dots + p^{rk}
(1+ \dots + p^{s-1})  $$
where $r$ is the largest integer less than $n/k$ and $s = n - rk$.  As $q$ divides every
term except possibly the last, $q$ divides $p^n-1$ if and only if it divides the last
term, that  is, if and only if $k$ divides $n$.  By taking $n$ prime, we ensure that no
prime $q$ which divides $1 + p + \dots + p^{k-1}$ for $k<n$ also divides $1 + p +
\dots + p^{n-1}$.  Thus $S$ is infinite.
\end{proof}

A minor modification to Example \ref{notpn} provides further examples of $p$-normal
spaces.

\begin{example}  Let $F$ be the field of $p^n$ elements.  Let $q$ be a prime dividing $p -
1$.  As before, we have $\xi \in F$ a primitive $q$-th root of unity and take $u$ a root
of $f(t) = t^q - c$ for some $c \in F$ which has no $q$-th root in $F$.  As before, the
space $\Lambda$ spanned by the conjugates of $u$ is $\langle u \rangle$.  We have $u^q =
c$.  We have $u^{p-1} = (u^q)^k = c^k$, where $k=(p-1)/q$, and $u^p = c^k u \in
\Lambda$.  It follows that $\Lambda^p = \Lambda$. \end{example}

\begin{definition} Let $\Lambda$ be a $p$-normal subspace of $\bar{F}$.  The
\textit{eigenvalue defined $p$-formation} $\pEv(\Lambda)$ is the class of all restricted
soluble Lie algebras $(L, \p)$ such that every eigenvalue of $\ad(x)$ is in $\Lambda$
for all $x \in L$. \end{definition}

The class $\pEv(\Lambda)$ is clearly a $p$-formation.  We shall show that it is
saturated.  For this, we need the following lemmas.

\begin{lemma} \label{allev} Let $L \in \Ev(\Lambda)$ be a soluble Lie algebra and let $V$
be an irreducible $L$-module giving the representation $\rho$.  Let $M$ be a maximal
ideal of $L$. Suppose all  eigenvalues of $\ad(b)$ are in $\Lambda$ for all $b \in M$,
but that for some $a \in L-M$, $\ad(a)$ has an eigenvalue not in $\Lambda$.  Then every
eigenvalue of $\ad(a)$ is outside $\Lambda$. \end{lemma}

\begin{proof} The characteristic polynomial $\chi(t)$ of $\rho(a)$ can be expressed as a
product $\chi(t) = f(t)g(t)$ where every root of $g(t)$ is in $\Lambda$ while no root of
$f(t)$ is in $\Lambda$.  Put $V_f = \{v \in V \mid f(\rho(a)) v = 0 \}$ and $V_g = \{v
\in V \mid g(\rho(a)) v = 0 \}$.  Then $V = V_f \oplus V_g$.  We want to prove that $V_g
= 0$.  As $V$ is irreducible, this follows if we prove that $V_g$ is an $L$-submodule of
$V$.  Consider $L$ and $V$ as $\langle a \rangle$-modules.  We work over the algebraic
closure, with the algebra $\bar{L} = \bar{F} \otimes L$ and module $\bar{V} = \bar{F}
\otimes V$.  Then $\bar{V}_g = \bar{F} \otimes V_g$.  

Let $\lambda$ be an eigenvalue of $\rho(a)$.  The weight space $\bar{V}_\lambda$ is the
space 
$$\{v \in \bar{V} \mid (\rho(a) - \lambda)^n v = 0 \text{ for some $n$}\}$$
and $\bar{V}_g = \sum_{\lambda \in \Lambda} \bar{V}_\lambda$.  Likewise, $\bar{L}$ is a
sum of weight spaces.  If $b \in \bar{L}_\mu$ for the eigenvalue $\mu$ of $\ad(a)$, and
$v \in \bar{V}_\lambda$, then $bv \in \bar{V}_{\lambda + \mu}$.  For $\lambda \in
\Lambda$, $\lambda + \mu \in \Lambda$ since all eigenvalues $\mu$ of $\ad(a)$ are in
$\Lambda$.  It follows that $b \bar{V}_g \subseteq \bar{V}_g$ for all $b \in \bar{L}$
and so, that $V_g$ is a submodule of $V$.  
\end{proof}

\begin{lemma} \label{evH} Let $(L, \p) \in \pEv(\Lambda)$ and let $V$ be an irreducible 
$(L, \p)$-module giving the representation $\rho$.  Suppose for some $x \in L$, $\rho(x)$
has an eigenvalue not in $\Lambda$.  Then $H^n(L, V) = 0$ for all $n$. \end{lemma}

\begin{proof}  Since any $L$-submodule of $V$ is an $(L, \p)$-submodule, $V$ is
irreducible as $L$-module.  We forget the $p$-operation and prove the result for
ordinary Lie algebras.  We use induction over $\dim(L)$.  The result holds if $\dim(L) =
1$.  Let $M$ be a maximal ideal of $L$ and let $W$ be a composition factor of $V$ as
$M$-module.  If there exists an element $b$ of $M$ for which $\rho(b)$ has an eigenvalue
outside $\Lambda$, then $\rho(b)| W$ has an eigenvalue outside $\Lambda$ since all
$M$-composition factors of $V$ are isomorphic.  By induction, $H^n(M,W) = 0$ for all $n$
and all composition factors $W$.  Therefore $H^n(M, V) = 0$ for all $n$.  By the
Hochschild-Serre spectral sequence, $H^n(L, V) = 0$.

We may therefore suppose that every eigenvalue of $\rho(b)$ is in $\Lambda$ for all $b
\in M$.  By Lemma \ref{allev}, we have $a \in L-M$ with every eigenvalue $\mu$ of
$\ad(a)$ outside $\Lambda$.  Every eigenvalue of the action of $a$ on the degree $q$
component $E^q(M)$ of the exterior algebra on $M$ is in $\Lambda$.  Thus every
eigenvalue of its action on $\Hom(E^q(M), V)$ is outside $\Lambda$.  As $H^q(M, V)$ is a
subquotient of $\Hom(E^q(M), V)$, every eigenvalue of the action of $a$ on $H^q(M, V)$ is
outside $\Lambda$.  Hence $H^p(L/M, W) = 0$ for every composition factor $W$ of
$H^q(M,V)$.  Thus $H^p(L/M, H^q(M, V)) = 0$ for all $p,q$ and the result follows by the
Hochschild-Serre spectral sequence. \end{proof}

\begin{theorem} \label{evdef} Let $\Lambda$ be a $p$-normal subspace of the algebraic
closure $\bar{F}$ of $F$.  Then the eigenvalue defined $p$-formation $\pEv(\Lambda)$ is
saturated. \end{theorem}

\begin{proof} Suppose $A$ is a minimal \p-ideal of $(L, \p)$, $(L/A, \p) \in
\pEv(\Lambda)$ but $(L, \p) \not\in \pEv(\Lambda)$.  We have to prove that $L$ splits
over $A$.  But this holds by Lemma \ref{evH}. \end{proof}

In characteristic $0$, every saturated formation is an eigenvalue defined formation.  If
we restrict attention to completely soluble algebras, this also holds in characteristic
$p \ne 0$ and for restricted algebras.

Let $\lambda \in \bar{F}$ have minimal polynomial $m(t)$ over $F$.  There exists a
vector space $V$ and linear transformation $a: V \to V$ with characteristic
polynomial $m(t)$.  Let $A$ be the Lie subalgebra $\langle a, a^p, a^{p^2}, \dots
\rangle$ of $\Hom(V, V)$.  Then the split extension of $V$ by $A$ is a primitive
completely soluble restricted Lie algebra which we denote by $P_\lambda$.  It is the
smallest algebra for which $\lambda$ appears as an eigenvalue.  (For the case of
ordinary Lie algebras, we use the ordinary primitive
algebra $\oP_\lambda$ the split extension of $V$ by $A = \langle a \rangle$.)  

\begin{lemma} \label{EvP}  Let $\f{F}$ be a saturated $p$-formation of completely soluble
algebras which contains all atoms.  Let $(L, \p) \in \f{F}$ and suppose the element $a
\in L$ has $\lambda$ as an eigenvalue of $\ad(a)$.  Then $P_\lambda \in \f{F}$. \end{lemma}

\begin{proof}  The element $a$ has $\lambda$ as an eigenvalue on some \p-chief factor
$V/W$ of $(L, \p)$.  As $V/W$ is $\f{F}$-central, the split extension of $V/W$ by
$L/\cent_L(V/W)$ is in $\f{F}$.  As $L$ is completely soluble, $L' \subseteq
\cent_L(V/W)$.  We may thus suppose $W = 0$ and that $L$ is the split extension of $V$
by an abelian algebra $M$ and that $M$ is a faithful irreducible $\f{F}$-central
$M$-module.  Let $\rho$ be the representation of $M$ on $V$.  We have $a \in M$ for
which $\rho(a)$ has $\lambda$ as an eigenvalue.  Note that every eigenvalue of $\rho(a)$
is a conjugate of $\lambda$.  Let $A = \langle a, a^\scp, a^{\scp^2}, \dots \rangle$.  We
construct the subdirect sum $M^*$ of two copies of $M$ by setting
$$M^* = \{(x,y) \in M \oplus M \mid (x + A) = (y + A) \}.$$
It is a restricted algebra with the $p$-operation $(x,y)^\scp = (x^\scp, y^\scp)$. 
We have projections $\pi_i: M^* \to M$ given by $\pi_1(x,y) = x$ and $\pi_2(x,y) =
y$.  Then $\ker \pi_1 = A_1 = \{(0, y) \mid y \in A\}$ and $\ker \pi_2 = A_2 = \{(x,0)
\mid x \in  \}$.  Since $M^*/A_i \simeq M$ and $A_1 \cap A_2 = 0$, $M^* \in \f{F}$.  The
diagonal subalgebra $D = \{(x,x) \mid x \in M\}$ is a \p-ideal since $M^*$ is abelian,
and we have $M^* = A_i \oplus D$.

We construct $M^*$-modules $V_1, V_2$ from copies of $V$ with action defined via $\pi_1,
\pi_2$, that is, $(x,y)v_1 = xv_1$ and $(x,y)v_2 = yv_2$ for $(x,y) \in M^*$ and $v_1 \in
V_1$, $v_2 \in V_2$.  Then $\cent_{M^*}(V_i) = A_i$ and $V_i$ is $\f{F}$-central.  The
$M^*$-module $W = \Hom(V_1, V_2)$ is $\f{F}$-hypercentral by Theorem \ref{tensHom}.  If
$f: V_1 \to V_2$ is an $M$-module homomorphism, then $((x,x)f)(v_1) = x \cdot f(v_1) -
f(xv_1) = 0$.  Thus $(x,x)f = 0$ for all $(x,x) \in D$ and $W^D \ne 0$.  But $W_D$ is
an $A_1$-module and is $\f{F}$-hypercentral.  Take any irreducible $A_1$-submodule $K$
of $W^D$ and form the split extension $P$ of $K$ by $A_1$.  Then $P \in \f{F}$.  Every
eigenvalue of $(a,0)$ on $V_1$ is a conjugate of $\lambda$, while $(a,0)V_2 = 0$.  Thus
every eigenvalue of $(a,0)$ on $K$ is a conjugate of $\lambda$.  As $K$ is irreducible
under the action of $a$, $P \simeq P_\lambda$.
\end{proof}

\begin{lemma} \label{allP}  Let $\f{F}$ be a saturated $p$-formation of completely soluble
algebras which contains all atoms.  Let $(A, \p)$ be abelian and let $V$ be a faithful
irreducible $(A, \p)$-module giving the representation $\rho$.  Let $P$ be the split
extension of $V$ by $A$.  Suppose that for all $a \in A$ and every eigenvalue $\lambda$
of $\rho(a)$, we have $P_\lambda \in \f{F}$.  Then $P \in \f{F}$. \end{lemma}

\begin{proof}  Take a basis $a_1, \dots, a_n$ of $A$.  Put $A_i = \langle a_i, a_i^\scp,
a_i^{\p^2}, \dots \rangle$.  Let $\lambda_i$ be an eigenvalue of $\rho(a_i)$ and let
$W_i$ be an $A_i$-module isomorphic to an irreducible $A_i$-submodule of $V$.  The split
extension of $W_i$ by $A_i$ is isomorphic to $P_{\lambda_i}$ and so, by assumption, is
in $\f{F}$.  Put $A^* = \oplus_i A_i$.  Then $W_i$ is an $A^*$-module with the summand
$A_j$ for $j \ne i$ acting trivially.  It is $\f{F}$-central.  We put $W = \otimes_i
W_i$.  Then $W$ is an $\f{F}$-hypercentral $A^*$-module by Theorem \ref{tensHom}.

We have a homomorphism $\phi: A^* \to A$ given by $\phi(b_1, \dots, b_n) = b_1 + \dots +
b_n$ for $b_i \in A_i$.  Now choose $v \in V$, $v \ne 0$ and $w_i \in W_i$, $w_i \ne
0$.  Put $\rho_i = \rho(a_i)$ and let $m_i(t)$ be the minimal polynomial of $\rho_i$. 
Any element of $W_i$ can be expressed in the form $f_i(\rho_i) w_i$ for some polynomial
$f_i(t)$ determined up to multiples of $m_i(t)$.  We have a map $\psi : W \to V$ defined
by
$$\psi(f_1(\rho_1)w_1 \otimes \dots \otimes f_n(\rho_n)w_n) =
f_1(\rho_1)  \dots f_n(\rho_n) v.$$
This is independent of the choices of the $f_i(t)$ as the $\rho_i$ commute and
$m_i(\rho_i) v = 0$.  Now $V$ is also an $A^*$-module via $\phi$ and $\psi$ is an
$A^*$-module homomorphism.  As $\psi$ is surjective, $V$ also is $\f{F}$-hypercentral
as $A^*$-module.  It is irreducible, so $\f{F}$-central and the split extension of $V$
by $A^*/ \cent_{A^*}(V)$ is in $\f{F}$.  But $\cent_{A^*}(V) = \ker\phi$.  Thus the split
extension is $P$.
\end{proof}

\begin{theorem} \label{everyEdef}  Let $\f{F}$ be a saturated $p$-formation of completely
soluble restricted Lie algebras which contains all atoms.  Then $\f{F} = \pEv(\Lambda)
\cap \pCs$ for some $p$-normal subspace $\Lambda$ of $\bar{F}$. \end{theorem}

\begin{proof}  Let $\lambda$ be the set of all eigenvalues of all elements of all
algebras $(L, \p) \in \f{F}$.  Then clearly $\f{F} \subseteq \pEv(\Lambda) \cap \pCs$. 
Suppose $(L,\p) \in \pEv(\Lambda) \cap \pCs$.  If $\lambda$ is an eigenvalue of $\ad(x)$
for some $x \in L$, then $\lambda \in \Lambda$ and by the definition of $\Lambda$,
is an eigenvalue of $\ad(A)$ for some element of an algebra in $\f{F}$.  By Lemma
\ref{EvP}, $P_\lambda \in \f{F}$.  By Lemma \ref{allP}, every primitive quotient of $(L,
\p)$ is in $\f{F}$.  Therefore $(L, \p) \in \f{F}$.
\end{proof}

\section{Intravariance of projectors}\label{intrav}
\begin{definition}  The \p-subalgebra $U$ of the restricted Lie algebra $(K, \p)$ is said to be
{\em intravariant} in $(K, \p)$ if every derivation of $K$ is expressible as the sum of an inner
derivation and a derivation which stabilises $U$.  \end{definition}

The $p$-operation plays no part in this definition, so $U$ is intravariant in $(K, \p)$ if and only if it is
intravariant in the underlying Lie algebra $K$.  By \cite[Lemma 1.2]{Frat}, the intravariant
\p-subalgebras of $(K, \p)$ are precisely those \p-subalgebras $U$ with the property that, if $K$ is a
\p-ideal of $(L, \p)$, then $L = K + \sN_L(U)$.  That a Cartan subalgebra $U$ of $(K, \p)$ (that is, a
nilpotent \p-subalgebra with $\sN_K(U) = U$) is intravariant is immediate from \cite[Theorem
2.1]{Frat} without any requirement that $K$ be soluble.  The analogue of \cite[Theorem 2.2]{Frat}
also holds with minor modification to cope with the existence of non-null atoms.

\begin{theorem} \label{thm-intrav}  Let $\f{H}$ be a $p$-homomorph of soluble restricted Lie
algebras which contains all atoms.  Let $K$ be a soluble \p-ideal of the restricted Lie algebra $(L, \p)$
and    let $S$ be an $\f{H}$-covering subalgebra of $K$.  Then $L = K + \sN_L(S)$. \end{theorem}

\begin{proof}  The result is trivial if $K = L$ or $S = K$.  Let $(L, \p)$ be a counterexample of least
possible dimension.  Let $A$ be a minimal \p-ideal of $(L, \p)$ contained in $K$.  Then $S+A/A$ is an
$\f{H}$-covering subalgebra of $K/A$.  By induction, we have
$$K/A + \sN_{L/A}(S+A/A) = L/A,$$
that is, $K + \sN_L(S+A) = L$.  Put $N = \sN_L(S+A)$.  Then $K \cap N$ is a \p-ideal of $N$ and $S$
is an $\f{H}$-covering subalgebra of $K \cap N$.  If $N < L$, then by induction we have $(K \cap N) + \sN_N(S)
= N$.  But then 
$$K + \sN_L(S) \supseteq K + (K \cap N) + \sN_U(S) = K + N = L$$
contrary to $(L, \p)$ being a counterexample.  Therefore $K_1 = S+A$ is a \p-ideal of $(L, \p)$.  It
is clearly sufficient to prove $K_1 + \sN_L(S) = L$.  As $K_1$ satisfies the conditions required of $K$,
we may replace $K$ with $K_1$, so we may suppose $S+A = K$.  We then have $K/ A \simeq S/S
\cap A \in \f{H}$.

By Corollary \ref{closure}, we may suppose that $\f{H}$ is saturated.  If $A \subseteq \Psi(L, \p)$,
then $K \in \f{H}$ by Theorem \ref{b-n} and $S =K$.  Therefore there exists a maximal
\p-subalgebra $U$ which complements $A$ in $L$.  Put $V = K \cap U$, $B = A \cap S$ and $T = B+
V$.
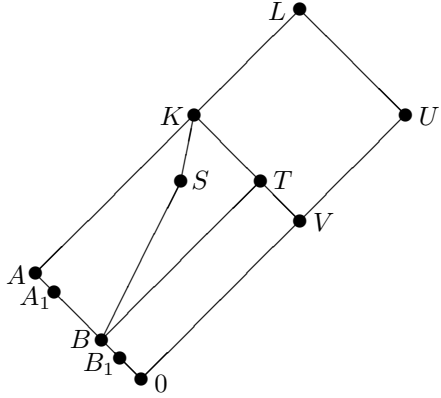
\begin{center}
\setlength{\unitlength}{1em}
\begin{picture}(16,16)(-7,0)
\put(3,15){\circle*{.5}}
\put(1.8,14.6){$L$}
\put(3,15){\line(-1,-1){10}}
\put(3,15){\line(1,-1){4}}
\put(7,11){\circle*{.5}}
\put(7.5,10.6){$U$}
\put(7,11){\line(-1,-1){10}}
\put(-3,1){\circle*{.5}}
\put(-2.5,0.5){$0$}
\put(-7,5){\line(1,-1){4}}
\put(-7,5){\circle*{.5}}
\put(-8.1,4.6){$A$}
\put(-1,11){\circle*{.5}}
\put(-2.3,10.6){$K$}
\put(-1,11){\line(1,-1){4}}
\put(3,7){\circle*{.5}}
\put(3.5,6.6){$V$}
\put(1.5,8.5){\circle*{.5}}
\put(2,8.2){$T$}
\put(-4.5,2.5){\line(1,1){6}}
\put(-4.5,2.5){\circle*{.5}}
\put(-5.7,2.2){$B$}
\put(-4.5,2.5){\line(1,2){3}}
\put(-1.5,8.5){\line(1,5){0.5}}
\put(-1.5,8.5){\circle*{.5}}
\put(-1.1,8.2){$S$}
\put(-3.8,1.8){\circle*{.5}}
\put(-5.2,1.3){$B_1$}
\put(-6.3,4.3){\circle*{.5}}
\put(-7.6,3.8){$A_1$}
\end{picture}
\end{center}

Since $A$ is abelian and $A+S = K$, $B$ is a \p-ideal of $K$.  Both $S$ and $T$ complement $A/B$
in $K/B$.  Since $S/B$ is an $\f{H}$-covering subalgebra of $K/B$, so is $T/B$.  By Lemma \ref{conj}, $T =
\alpha_a(S)$ for some $a \in A$.  But $\alpha_a = 1 + \ad(a)$ is an automorphism of $(L, \p)$, so we
may suppose $S = T = B+V$.

If $B=0$, then $S=V$, $\sN_L(S) \supseteq U$ and $K+ \sN_L(S) = L$.  Therefore $B \ne 0$.  Let
$B_1$ be a maximal $K$-submodule of $B$.  Then $S/B_1$ is the split extension of $B/B_1$ by $V$
and is in $\f{H}$.  Let $A_1$ be a maximal $K$-submodule of $A$ containing $B$.  Since $A$ is
irreducible as $L$-module and $K$ is an ideal of $L$, all $K$-composition factors of $A$ are
isomorphic by Zassenhaus \cite[Lemma 1]{Zas}.  Since both $K/A_1$ and $S/B_1$ are split
extensions of composition factors of $A$ by $V$, they are isomorphic and $K/A_1 \in \f{H}$.  But $S
+ A \ne K$ contrary to $S$ being an $\f{H}$-covering subalgebra of $K$.  Thus no minimal counterexample
exists.  \end{proof}

\section{Comparisons}\label{compar}
In this section, I investigate the relationship between Schunck classes of soluble Lie algebras  and
$p$-Schunck classes of soluble restricted Lie algebras over the same field $F$ of characteristic $p$.  
For this, I define a function $\Res$ classes of Lie algebras to classes of restricted Lie algebras and a
function $\Und$ from classes of restricted Lie algebras to classes of Lie algebras.

\begin{definition} Let $\f{H}$ be a class of soluble Lie algebras over
$F$.  We define $\Res(\f{H})$ to be the class of all restricted Lie algebras $(L, \p)$
with $L \in \f{H}$. For a class $\f{K}$ of soluble restricted Lie algebras, we define
$\Und(\f{K})$ to be the class of underlying algebras of members of $\f{K}$.
\end{definition}

\begin{lemma} \label{lem-res} Let $\f{H}$ be a homomorph.  Then $\Res(\f{H})$ is a
$p$-homomorph.  If $\f{H}$ is a formation, then $\Res(\f{H})$ is a $p$-formation. 
\end{lemma}

\begin{proof} If $(L, \p) \in \f{H}$ and $A$ is a \p-ideal of $L$, then $L/A \in \f{H}$,
so $(L/A, \p) \in \Res(\f{H})$.  If $\f{H}$ is a formation and $(L, \p)$ is a restricted
Lie algebra with $p$-ideals $A, B$ such that $(L/A, \p)$ and  $(L/B, \p)$ are in $\Res(\f{H})$,
then $L/A$ and $L/B$ are in $\f{H}$. So $L/(A \cap B) \in \f{H}$ and $(L/(A
\cap B), \p) \in \Res(\f{H})$. \end{proof}

Note that in the above, we had a restricted Lie algebra  $(L, \p)$ and quotients $L/A,
L/B$ in the formation $\f{H}$.  For a Lie algebra $L$, having restrictable quotients
$L/A, L/B$ does not imply $L/(A \cap B)$ restrictable.

\begin{example}\label{noform} Let $U = \langle a,b,c \rangle$ with $ab=b. ac=bc=0$.  Let $V =
\langle v_0, \dots, v_{p-1}$ be the $U$-module with the action $av_i = iv_i, bv_i = v_{i+1}$, where
the indices are integers mod $p$.  Let $W = \langle w_0, w_1\rangle$ be the $U$-module with action
$aw_0 = 0, aw_1 = w_1, bw_0 = -w_1, bw_1 = 0, cw_i = w_i$.  Let $X$ be the split extension of $V$ by
$U$.  Putting $a^\scp = a, b^\scp = c^\scp = c, v_i^\scp =0$ makes $(X, \p)$ a restricted Lie algebra.  Let $Y$
be the split extension of $W$ by $U$.  Putting $a^\scpd = a, b^\scpd = 0, c^\scpd = c, w_i^\scpd = 
0$ makes $(Y, \pd)$ a restricted Lie algebra.  Since $X$ and $Y$ have trivial centres, \p and \pd are the
only $p$-operations on them.  Let $L$ be the split extension of $V \oplus W$ by $U$.  Any $p$-operation
on $L$ would have to agree with $\p$ on $L/W$ and with $\pd$ on $L/V$.  But $\p$ and $\pd$ do not
agree on $U = L/(V+W)$.
\end{example}

From Example \ref{noform}, we see that $\f{K}$ being a $p$-formation which contains with any
$(L,\p)$, also $(L,\pd)$ for any $p$-operation $\pd$ on $L$ does not ensure that $\Und(\f{K})$ is a
formation. 

\begin{lemma}\label{lem-und} Let $\f{K}$ be a $p$-homomorph with the property that, if $(L,\p) \in
\f{K}$ and \pd is another $p$-operation on $L$, then $(L,\d) \in \f{K}$.  Then $\Und(\f{K})$ is a
homomorph. \end{lemma}

\begin{proof} It is sufficient to prove that if $A$ is a minimal ideal of $L \in \Und(\f{K})$, then $L/A \in
\Und(\f{K})$.  By Lemma \ref{abpdid}, $A$ is a \pd-ideal  for some $p$-operation \pd on $L$, 
so $L/A \in \Und(\f{K})$. \end{proof}

However, $\f{K}$ a $p$-Schunck class does not ensure that $\Und(\f{K})$ is a Schunck class.  If the
$p$-Schunck class $\f{K}$ contains the algebras $(X, \p)$ and $(Y, \pd)$ of Example \ref{noform}, then
$\Und(\f{K})$ contains  $X$ and $Y$.  If it is a Schunck class, then it contains any algebra whose
primitive quotients are quotients of $X$ or $Y$.  In particular, it must contain $L$ which, being not
restrictable, cannot be in $\Und(\f{K})$.  I introduce another function which does give Schunck classes.

\begin{definition} Let $\f{K}$ be a $p$-Schunck class.  We define $\Ord(\f{K}) = \prim(\Und(\f{K}))$.
\end{definition}

As every non-zero Schunck class contains the class $\f{N}$ of nilpotent Lie algebras, and every
$p$-Schunck class containing all atoms contains the class $p\f{N}$, we restrict attention to
($p$-)Schunck classes containing $(p)\f{N}$.  We denote the set of all  non-zero Schunck classes by
$\Sch$ and the set of all $p$-Schunck classes which contain all atoms by $\pSch$.  They are partially
ordered by inclusion.  We make them into lattices by defining the sum of ($p$-)Schunck classes.

\begin{definition} Let $\f{H}$ be a ($p$-)homomorph.  We define the {\em skeleton}
$\Skel(\f{H})$ of $\f{H}$ to be the class of all primitive algebras in $\f{H}$.   We
define the {\em sum} or {\em join} of two ($p$-)Schunck  classes $\f{H}_1, \f{H}_2$
by
$$ \f{H}_1 + \f{H}_2 = \prim(\Skel(\f{H}_1) \cup \Skel(\f{H}_2)).$$
\end{definition} 
For a ($p$-)Schunck class $\f{H}$, we clearly have $\prim(\Skel(\f{H})) = \f{H}$ and $\f{H}_1 +
\f{H}_2$ is the smallest ($p$-)Schunck class containing both $\f{H}_1$ and $\f{H}_2$.  The 
($p$-)Schunck classes are in one-one correspondence with their skeleta and the operations $\cap, +$
on them correspond to the set-theoretic operations $\cap, \cup$ on the skeleta.  Thus \Sch\ and    
\pSch\ are lattices.

\begin{lemma}  \Sch\ and \pSch, partially ordered by inclusion, are complete
distributive lattices. \end{lemma}

\begin{proof} If $L$ is in some infinite sum of Schunck classes, it has
only finitely many chief factors, so only finitely many isomorphism types of primitive
quotients.  It thus is in the sum of a finite subset of the Schunck classes.  Thus
$$ \sum_i \f{H}_i = \prim\bigl(\bigcup_i \Skel(\f{H}_i)\bigr),$$
and the result for \Sch\ follows. The result for \pSch\ follows similarly.\end{proof}

Note that \Sch\ and \pSch\ have greatest (the classes $\f{S}, p\f{S}$ of all soluble
algebras) and least elements (the classes $\f{N}, p\f{N}$ of all nilpotent algebras).

The intersection of two saturated $p$-formations is a saturated $p$-formation, however their
sum need not be a $p$-formation.

\begin{example}  Let $\Lambda_1, \Lambda_2$ be $p$-normal subspaces of $\bar{F}$,
neither of which contains the other (for example, normal extension fields of $F$ of
relatively prime degrees).  There exist $\lambda_i \in \Lambda_i$ for which $\lambda_1 +
\lambda_2$ is in neither space.  We have the primitive algebras $P_{\lambda_i} \in
\pEv(\Lambda_i)$.  Any saturated $p$-formation which contains both $P_{\lambda_1}$ and
$P_{\lambda_2}$ must also contain $P_{\lambda_1+\lambda_2}$.  But 
$P_{\lambda_1+\lambda_2} \not\in \pEv(\Lambda_1) + \pEv(\Lambda_2)$. \end{example}

The above construction is not possible over an algebraically closed field.  To provide an
example over an arbitrary field, I use the standard examples of soluble algebras with
non-nilpotent derived algebras.  Let $N = \langle a,b,c \rangle$ be the nilpotent Lie
algebra, $ab = c$, $ac = bc = 0$ and let $K = \langle k_0, \dots, k_{p-1} \rangle$, where
the indices are integers mod $p$, be the $N$-module with $a k_i = i k_{i-1}$, $b k_i =
k_{i+1}$ and $c k_i = k_i$.  We set $a^\scp = 0$, $b^\scp = c$ and $c^\scp = c$.  Then $K$
is an $(N,\p)$-module.  Let $P$ be the split extension of $K$ by $N$.  Let $S =
\langle x, y \rangle$ be the Lie algebra with $xy = y$.  Let $V = \langle v_0, \dots,
v_{p-1} \rangle$ be the $S$-module with $x v_i = i v_i$ and $y v_i = v_{i+1}$ and let
$Q$ be the split extension of $V$ by $S$.  Then $Q$ is not restrictable, so we also
use the algebra $Q^*$ which is the split extension of $V$ by $S^* = S \oplus \langle z
\rangle$ with $z v_i = v_i$.

\begin{example}\label{explus}  Let $\f{F}_1$ be the saturated $p$-formation generated by $P$
and let
$\f{F}_2$ be the saturated $p$-formation generated by $Q^*$.  Suppose $\f{F}$ is a
saturated $p$-formation which contains both $P$ and $Q^*$.  Then $N \oplus S^* \in
\f{F}$.  With $S^*$ acting trivially on $K$ and $N$ acting trivially on $V$, 
$K$ and $V$ are $\f{F}$-central $(N \oplus S^*)$-modules.  By Theorem \ref{tensHom}, $K
\otimes V$ is $\f{F}$-hypercentral.  It is a faithful irreducible $(N \oplus
S^*)$-module, so the split extension $T$ is a primitive algebra in $\f{F}$.  Denoting
the $p$-formation of algebras which are nilpotent of class at most 2 by $p\f{N}_2$, we
have $P \in \f{F}^*_1 = \pLoc(p\f{N}_2)$, so $\f{F}_1 \subseteq \f{F}^*_1$.  Let
$\f{M}$ be the $p$-formation of metabelian algebras with the property that every
nilpotent subalgebra is abelian.  Then $Q^* \in \f{F}^*_2 = \pLoc(\f{M})$ and
$\f{F}_2 \subseteq \f{F}^*_2$.  But $T$ is primitive and not in either $\f{F}^*_i$, so $T
\not\in \f{F}_1 + \f{F}_2$. \end{example}  

\begin{lemma} \label{Res}   Let $\f{H} \ne 0$ be a Schunck class of soluble Lie
algebras  over $F$.  Then $\Res(\f{H})$ is a $p$-Schunck class which contains every
atom.  \end{lemma} 

\begin{proof}   Every abelian algebra is in $\f{H}$, so every atom is in $\Res(\f{H})$. 
Suppose every primitive quotient of $(L,\p)$ is in $\Res(\f{H})$.  Let $K$ be an ideal of $L$ with
$L/K$ primitive.  If $L/K$ is abelian, then $L/K \in \f{H}$.  If $L/K$ is non-abelian, then by Lemma
\ref{ordprim}, $K$ is a \p-ideal of $L$, so $(L/K,\p) \in \Res(\f{H})$ and $L/K \in \f{H}$.  Thus $L \in
\prim\f{H} =  \f{H}$ and $(L,\p) \in\Res(\f{H})$.  
\end{proof}

\begin{theorem} \label{lattice}  $\Res: \Sch \to \pSch$ and $\Ord: \pSch \to
\Sch$ are homomorphisms of complete lattices and, for any $p$-Schunck class $\f{K}$,
$$\Res(\Ord(\f{K})) = \f{K}.$$ \end{theorem}

\begin{proof}  The lattice operations correspond to set operations on the skeleta and
the result follows. \end{proof}

\section{$p$-envelopes}\label{env}
We use the basic theory of $p$-envelopes as set out in Strade and Farnsteiner \cite[pp. 94--97]{SF}.

\begin{lemma} \label{env-null}  Let $A$ be an abelian ideal of the soluble Lie algebra $U$.  Then there
exists a $p$-envelope $(L, \p)$ of $U$ in which $A$ is a null \p-ideal. \end{lemma}

\begin{proof} Take any $p$-envelope $(L, \p)$ of $U$.  By Lemma \ref{abpdid}, there exists a
$p$-operation   
\pd on $L$ with $A^\scpd = 0$.    Replacing $L$ by the \pd-closure $U_\scpd$ of $U$ gives the required
$p$-envelope. \end{proof}

\begin{lemma} \label{prim1} Suppose $P$ is a primitive soluble Lie algebra and let $(L, \p)$ be a minimal
$p$-envelope of $P$.  Then $(L,\p)$ is primitive. \end{lemma}

\begin{proof} If $P$ is abelian, then it is 1-dimensional and so is any minimal $p$-envelope.  Suppose
$P$ is non-abelian.  By Strade and Farnsteiner \cite[Theorem 5.8, p. 96]{SF}, $\z(L) \subseteq P$.  But
$\z(P) = 0$, so $\z(L) = 0$.  Let $A$ be the minimal ideal of $P$.  Then $\ad(a)^2 = 0$ for all $a \in
A$, so $A$ is a null \p-ideal of $(L, \p)$.  We have $H^\beta(P/A, A) = 0$ for all $\beta$, so $H^n(L/A,
A) = 0$ for all $n$.  Thus there exists a complement $M$ to $A$ in $L$.  Suppose $M$ contains a minimal
ideal $C$ of $L$.  Then $C \cap P = 0$ as $A$ is the only minimal ideal of $P$.  Therefore $C \subseteq
\z(L)$ by \cite[Lemma 5.5, p. 94]{SF}.  But $\z(L) = 0$.  Therefore $\cent_L(A) = A$. \end{proof}

Note that, by \cite[Theorem 5.8, p. 96]{SF}, any two minimal $p$-envelopes are isomorphic as Lie
algebras.  As the minimal $p$-envelopes of a non-abelian primitive algebra have trivial centre, they are
isomorphic as restricted Lie algebras.

\begin{lemma} \label{prim2} Suppose $(L, \p)$ is a $p$-envelope of the non-abelian Lie algebra     
$U$ and that $(L, \p)$ is primitive.  Then $U$ is primitive. \end{lemma}

\begin{proof} Let $A$ be a minimal ideal of $U$.  Then $A$ is an abelian ideal of $L$, so   $A^\scp
\subseteq \z(L) = 0$.  Therefore $A$ is the unique minimal \p-ideal of $L$.  Thus $\cent_U(A) = A$ and     
$U$ is primitive.\end{proof}

Note that it is possible for a restricted Lie algebra to be the minimal $p$-envelope of distinct
primitive algebras.  If $(L,\p)$ is a minimal $p$-envelope of a non-restrictable primitive Lie algebra
$U$, then it is also a minimal $p$-envelope of the primitve algebra $L$.

\begin{definition}   Let $\f{K}$ be a $p$-Schunck class of soluble restricted Lie algebras which contains all
atoms.  We define the {\em enveloped class} $\envd(\f{K})$  of $\f{K}$ to be the class of all Lie
algebras $L$ having a $p$-envelope in $\f{K}$. \end{definition}

\begin{lemma}\label{mine}  Suppose $\f{K}$ is a $p$-Schunck class which contains all atoms.  
Suppose $(L, \p) \in \f{K}$ is a $p$-envelope of $U$.  Then every minimal $p$-envelope of $U$ is in
$\f{K}$. \end{lemma}

\begin{proof}  Let $i: U \to L$ be the inclusion.  Let $(M,\pd)$  be a minimal $p$-envelope of $U$ with
$i': U \to M$ the inclusion.  By \cite[Proposition 5.6, p. 95]{SF}, there exists a homomorphism $f: L \to M$ with
$f \circ i = i'$.  Let $K = \ker(f)$.  Then $K$ is an ideal of $L$ and $K \cap U = 0$.  By \cite[Lemma 5.5,
p.94]{SF}, $K \subseteq \z(L)$.  Further, $L/K$ with some $p$-operation $\p''$ is a $p$-envelope of $U$.  But
$M$ has the least possible dimension for a $p$-envelope, so $f(L) = M$.  Now for some $p$-operation $\p^*$
on $L$, $K$ is a $\p^*$-ideal.  By Theorem \ref{anyop}, $(L, \p^*) \in \f{K}$, so $(L/K, \p^*) \in \f{K}$,
that is, $(M, \p^*) \in \f{K}$.  By Theorem \ref{anyop}, $(M, \pd) \in \f{K}$.  \end{proof}

\begin{lemma}\label{anye}  Suppose $\f{K}$ is a $p$-Schunck class which contains all atoms. 
Suppose $U$ has a $p$-envelope in $\f{K}$.  Then every $p$-envelope of $U$ is in $\f{K}$.
\end{lemma}

\begin{proof}  By Lemma \ref{mine}, $U$ has a minimal $p$-envelope $(M,\pd) \in \f{K}$.  Let $(L, \p, i)$
be a $p$-envelope.  There exists a homomorphism $f: L \to M$ with $f \circ i = i'$.  Let $K = \ker(f)$. 
Then $f(L) = M$.  $K$ is central in $L$ and is a $\p^*$-ideal for some $p$-operation $\p^*$.  We have
$(L/K, \p^*) \in \f{K}$ and so $(L, \p^*) \in \f{K}$ by Corollary \ref{central}.  Therefore also $(L,
\p) \in \f{K}$. \end{proof}

\begin{lemma}\label{envdH} Suppose $\f{K}$ is a $p$-Schunck class which contains all atoms. 
Then $\envd(\f{K})$ is a homomorph. \end{lemma}

\begin{proof}  Suppose $U \in \envd(\f{K})$.  Let $A$ be a minimal ideal of $U$.  By Lemmas
\ref{env-null} and \ref{anye}, we can choose a $p$-envelope $(L, \p) \in \f{K}$ with $A$ a \p-ideal. 
Then $(L/A, \p) \in \f{K}$ is a $p$-envelope of $U/A$, so $U/A \in \envd(\f{K})$.  It follows that $U/K
\in \f{K}$ for any ideal $K$ of $U$. \end{proof}

\begin{theorem}\label{envdSat} Suppose $\f{K}$ is a $p$-Schunck class which contains all
atoms.  Then $\envd(\f{K})$ is a Schunck class. \end{theorem}

\begin{proof} Suppose $A$ is a minimal ideal of $U$, $U/A \in \envd(\f{K})$ and $U \not\in
\envd(\f{K})$.  Take a minimal $p$-envelope $(L, \p)$ of $U$ such that $A$ is a \p-ideal.  Now $(L/A,
\p) \in \f{K}$ by Lemma \ref{anye}, but $(L, \p) \not\in \f{K}$.  There exists a $\f{K}$-projector $M$
of $(L, \p)$.  $M$ complements $A$ in $L$ and $M \cap U$ complements $A$ in $U$.  We have to show
that $M \cap U$ is an $\envd(\f{K})$-projector of $U$.

Let $K$ be an ideal of $U$.  We have to show that, if $U/K \in \envd(\f{K})$, then $(M \cap U) + K =
U$.  This holds if $K \not\subseteq M$, so we may suppose $K \subseteq M$.  Now $K$ is an ideal of
$L$.  Let $B \subseteq K$ be a minimal ideal of $L$.  There exists a $p$-operation \pd on $L$ which
vanishes on the abelian ideal $A+B$.  From the minimality of $(L, \p)$, $(L, \pd)$ is also a minimal
$p$-envelope of $U$, so we may assume that $B$ is a \p-ideal.   Since $U/K \in \envd(\f{K})$, we have
$(L/K, \p) \in \f{K}$ by Lemma \ref{anye}.   But $M$ is a $\f{K}$-projector of $L$ and $M+K = L$
contrary to $K \subseteq M$.
\end{proof}

\begin{theorem}\label{lhomom} The class map $\envd : \pSch \to \Sch$ is a lattice
homomorphism which sends $p$-formations to formations. \end{theorem}

\begin{proof} Let $\f{K}, \f{K}' \in \pSch$.  Suppose $U \in \envd(\f{K}) \cap \envd(\f{K}')$.  Then
$U$ has a $p$-envelope $(L,\p) \in \f{K}$ and a $p$-envelope $(L',\pd) \in \f{K}'$.  By Lemma
\ref{anye},  $(L,\p) \in \f{K}'$, so $U \in \envd(\f{K} \cap \f{K}')$.  Trivially, if $U \in \envd(\f{K} \cap
\f{K}')$, then $U \in \envd(\f{K}) \cap \envd(\f{K}')$.  Thus $\envd(\f{K} \cap \f{K}') = \envd(\f{K}
\cap \f{K}')$.

Suppose $U \in \envd(\f{K}) \cup\, \envd(\f{K}')$.  Let $(L, \p)$ be a $p$-envelope of $U$.  Let $(L/K,
\p)$ be a primitive quotient of $(L,\p)$.  Then $(L/K,\p)$ is a $p$-envelope of $U+K/K$.   If $U+K/K$
is abelian, $(L/K,\p)$ is an atom and by assumption, is in both $\f{K}$ and $\f{K}'$.  Suppose
$U+K/K$ is non-abelian.  Then  by Lemma \ref{prim2}, $U+K/K$ is primitive.  But every primitive
quotient of $U$ is in $\envd(\f{K})$ or $\envd(\f{K}')$ and it follows that $(L/K,\p)$ is in either
$\f{K}$ or $\f{K}'$.  Thus $(L,\p) \in \f{K} \cap \f{K}'$ and $U \in \envd(\f{K} \cap \f{K}')$.

Suppose $U \in \envd(\f{K} \cup \f{K}')$.  Let $K$ be an ideal of $U$ with $U/K$ primitive.  If $K = 0$,
then a minimal $p$-envelope of $U$ is primitive, is in $\envd(\f{K} \cup \f{K}')$ and so in $\f{K}$ or
$\f{K}'$.  We then have $U \in \envd(\f{K}) \cup \envd(\f{K}')$, so suppose $K \ne 0$.  Let $A
\subseteq K$ be a minimal ideal of $U$.  By Lemma \ref{env-null}, there exists a $p$-envelope
$(L,\p)$ of $U$ in which $A$ is a null \p-ideal.  It follows that $(L/A, \p) \in \f{K} \cup \f{K}'$ is a
$p$-envelope of $U/A$.  By induction on $\dim(U)$, $U/K$ is in $\envd(\f{K})$ or in $\envd(\f{K}')$.  
Thus $U \in \envd(\f{K}) \cup \envd(\f{K}')$. 

Now suppose $\f{K} \in \pSch$ is a $p$-formation.   Suppose $A_1, A_2$ are minimal ideals of $U$ and
that $U/A_i \in \envd(\f{K})$.  We can choose a $p$-envelope $(L,\p)$ of $U$ such that $\p$ vanishes
on the abelian ideal $A_1 + A_2$.  We then have $(L/A_i, \p)$ is a $p$-envelope of $U/A_i$ and so is in
$\f{K}$.  Therefore $(L, \p) \in \f{K}$ and $U \in \envd(\f{K})$.
\end{proof}


\begin{thebibliography}{00}

\bibitem{CoSol} D. W. Barnes, \textit{On the cohomology of soluble Lie algebras} 
Math. Zeitschr. \textbf{101} (1967), 343--349


\bibitem{Frat} D. W. Barnes,  \textit{The Frattini argument for Lie algebras,}  Math. Zeitschr. \textbf{133}
(1973), 277--283.


\bibitem{HyperC} D. W. Barnes,  \textit{On $\fF$-hypercentral modules for  Lie
algebras,}  Archiv der Math. \textbf{30} (1978), 1--7.


\bibitem{excentric} D. W. Barnes,  \textit{On $\mathfrak{F}$-hyperexcentric modules
for Lie algebras,}  J. Austral.  Math. Soc.  \textbf{74} (2003), 235--238.


\bibitem{extras} D. W. Barnes,  \textit{Ado-Iwasawa extras,}  J. Austral.  Math. Soc.  to appear.

\bibitem{BGH} D. W. Barnes and H. M. Gastineau-Hills,  \textit{On the theory of soluble Lie
algebras,}  Math. Zeitschr.  \textbf{106} (1968), 343--354.

\bibitem{BN} D. W. Barnes and M. L. Newell,  \textit{Some theorems on
saturated homomorphs of soluble Lie algebras,}  Math. Zeitschr. \textbf{115} (1970), 179--187.

\bibitem{DG} M. Demazure and P. Gabriel, \textit{Introduction to algebraic geometry and algebraic
groups,}  North-Holland Mathematics Studies \textbf{39}, North-Holland, Amsterdam-New
York-Oxford, 1980.

\bibitem{DH} K. Doerk and T. Hawkes, \textit{Finite soluble groups,} DeGruyter,  Berlin-New York
1992.

\bibitem{H} G. Hochschild, \textit{Cohomology of restricted Lie algebras,} Amer. J. Math. \textbf{76}
(1954), 555--580.

\bibitem{Jac} N. Jacobson, \textit{Lie algebras}  Interscience, New York-London, 1962.

\bibitem{Sch} E. Schenkman,  \textit{A theory of subinvariant Lie algebras,}  Amer. J.
Math.  \textbf{73} (1951), 453--474.

\bibitem{SF} H. Strade and R. Farnsteiner, \textit{Modular Lie algebras and their
representations},  Marcel  Dekker, Inc., New York-Basel, 1988.

\bibitem{Stit1} E. L. Stitzinger, \textit{Covering-avoidance for saturated formations of solvable Lie
algebras,} Math. Zeitschr. \textbf{124} (1972), 237--249.

\bibitem{Stit2} E. L. Stitzinger, \textit{On saturated formations of solvable Lie algebras,} Pacific J.
Math. \textbf{47} (1973), 531--538.

\bibitem{Vgt} D. Voigt, \textit{Induzierte Darstellungen in der Theorie der endlichen,
algebraischen Gruppen,} Lecture Notes in Mathematics \textbf{592} Springer-Verlag,                  
Berlin-Heidelberg-New York, 1977.

\bibitem{Zas} H. Zassenhaus, \textit{On trace bilinear forms on Lie-algebras,} Proc. Glasgow Math.
Assoc. \textbf{4} (1959), 62--72.

\end{thebibliography}
\end{document}